\documentclass[a4paper,12pt]{article}
\usepackage{authblk}
\usepackage[latin1]{inputenc}
\usepackage{amsfonts,amssymb,amsmath,amsthm}
\usepackage{pdfsync}
\usepackage{color}
\usepackage{url}
\usepackage{nicefrac}
\usepackage[pdftex]{hyperref}
\usepackage{hyperref}

\numberwithin{equation}{section}
\hypersetup{colorlinks=true,citecolor=black, filecolor=black, linkcolor=black, urlcolor=black}

\makeatletter
\@addtoreset{equation}{section}

\makeatother

\newtheorem{theo}{Theorem}[section]
\newtheorem{lemm}[theo]{Lemma}
\newtheorem{prop}[theo]{Proposition}
\newtheorem{rema}[theo]{Remark}
\newtheorem{coro}[theo]{Corollary}

\newcommand{\pA}{\partial A}
\newcommand{\eps}{\varepsilon}

\newcommand{\R}{\mathbb{R}}
\newcommand{\N}{\mathbb{N}}
\renewcommand{\L}{\mathcal{L}}
\newcommand{\E}{\mathcal{E}}
\newcommand{\D}{\mathcal{D}}
\renewcommand{\H}{\mathcal{H}}

\renewcommand{\O}{\mathcal{O}}

\newcommand{\F}{\mathcal{F}}

\newcommand{\G}{\Gamma}

\renewcommand{\l}{\ell}
\renewcommand{\a}{\alpha}

\renewcommand{\S}{\Sigma}

\newcommand{\vp}{\varphi}
\newcommand{\vpeps}{\varphi_\eps}

\newcommand{\pD}{\partial \D}
\newcommand{\mesrest}{\, \text{\Large$\llcorner$}\, }

\newcommand{\be}{\begin{equation}}
\newcommand{\ee}{\end{equation}}
\newcommand{\ben}{\begin{equation*}}
\newcommand{\een}{\end{equation*}}
\newcommand{\ba}{\begin{eqnarray}}
\newcommand{\ea}{\end{eqnarray}}
\newcommand{\ban}{\begin{eqnarray*}}
\newcommand{\ean}{\end{eqnarray*}}

\begin{document}

\title{A minimal interface problem arising from a two component Bose Einstein condensate via $\G$-convergence}
\author{Amandine Aftalion${}^{1}$ \& Jimena Royo-Letelier${}^{1,} {}^{2}$}
\affil{${}^1$
Université de Versailles Saint-Quentin, CNRS UMR 8100,\\Laboratoire de Math\'ematiques de Versailles, \\ 45
avenue des \'Etats-Unis, 78035 Versailles C\'edex, France \\
${}^2$Ceremade, CNRS UMR 7534, Université Paris-Dauphine, Place du Maréchal
de Lattre de Tassigny, 75775 Paris Cédex 16, France.}
\date{\today}

\maketitle

%Keywords: Phase transition, Two-component Bose Einstein condensate, $\Gamma$ convergence

\begin{abstract}
 We consider the energy modeling a two component
 Bose-Einstein condensate in the limit of strong coupling and strong segregation. We prove the $\Gamma$-convergence to a perimeter minimization problem, with a weight given by the density of the condensate. In the case of equal mass for
  the two components, this leads to symmetry breaking for the ground state. The proof relies on a new formulation of the problem in terms of the total density and spin functions, which turns the energy into
  the sum of two weighted Cahn-Hilliard energies. Then, we use techniques coming from geometric measure theory
    to construct upper and lower bounds. In particular, we make use of the slicing technique introduced in  \cite{AT}.
\end{abstract}

%%\tableofcontents

%%%%%%%%%%%%%%%%
%%%%%%%%%%%%%%%%
\section{Introduction} \label{intro}

The aim of this paper is to prove a $\Gamma$-convergence result for a functional modeling a two component
 Bose-Einstein condensate in the case of segregation. We introduce a new formulation of the problem which transforms the two wave functions describing each component
 of the condensate into total density and spin functions. The new functional in the density and spin variables is given by the sum of two weighted Cahn-Hilliard energies modeling phase transition problems as in the  Modica-Mortola problem \cite{Mod}. In fact, our new functional is strongly related  to that of Ambrosio-Tortorelli approaching the Mumford-Shah image segmentation functional \cite{AT}.   We use techniques coming from geometric measure theory \cite{alb,ABS,AT,BOU} to construct upper and lower bounds for our initial functional and prove $\Gamma$-convergence to a perimeter minimization problem, with a weight given by the density of the condensate. There is a large mathematical literature about the segregation patterns for two component Bose Einstein condensates \cite{BeLinWeiZhao,BeTer,CaffLin2,ctv3,NoTaTeVe,WeWe1}: regularity of the limiting functions, regularity of the interface, asymptotic behaviour near the interface. All these papers use the limiting equations and do not take into account
  the trapping potentials and the $\Gamma$ convergence of the energy as we do.

Before introducing the functional for a two component Bose Einstein condensate, we recall some properties of a single Bose Einstein condensate (BEC). A single BEC is described by the wave function $\eta$ minimizing the energy
  \be\label{enetaeps}
	 E_\eps(\eta) = \frac12 \int_{\R^2} |\nabla \eta|^2 + \frac1{\eps^2} V(x) |\eta|^2  +   \frac1{2\eps^2} |\eta|^4   \,
\ee
where $V$ is the trapping potential, usually taken to be harmonic, that is $V(x)=|x|^2$, $\eps$ is a small parameter
 giving rise to a large coupling constant describing the repulsive self interaction of the condensate. The minimization
 is performed under the mass constraint $ \int_{\R^2} |\eta |^2 = 1$. We define the  ground state by

 \be \label{definitionetaeps}
	E_\eps(\eta_\eps)  =  \inf_{ \int_{\R^2} |\eta |^2 = 1 }  E_\eps(\eta  ) \,,
	 \ee

which is, up to multiplication
  by a constant,  a real positive function. Let \be\label{rho}\rho(x) = \max ( \lambda^2 -  |x|^2,0) \hbox{ with }\lambda >0 \hbox{ chosen such that }\int_\D \rho=1\hbox{ where }\D = B(0,\lambda).\end{equation} Then, when $\eps$ is small, the ground state $\eta_\eps$
   is close to the function $\sqrt \rho$ in $\D$, with exponential decay at infinity. Properties of $\eta_\eps$ can be found in
   \cite{Aflivre,AJR,GaPe,IM,KaSou}.
   %When the condensate is put under rotation, singularities called vortices exist
   %above a critical value of the velocity. The velocity $\Omega_c$ where it is energetically favorable to have vortices has been
   %analyzed in \cite{IM}. The first vortex appears close to the origin and in \cite{AJR}, it is proved that no vortex
   %exists, even in the low density region $\R\backslash \D$, and that up to $\Omega_c$, the ground state is real valued.

   A two component Bose Einstein condensate can be experimentally realized as 2 isotopes of the same atom in different
   spin states \cite{hall} or isotopes of different atoms \cite{modugno}.
    They are described by two wave functions $u_1$ and $u_2$, respectively representing components 1 and 2. The Gross Pitaevskii
    energy  of the two component condensate is given by
\be\label{enertwobec}
	 \E_\eps(u_1, u_2)  = E_\eps(u_1) + E_\eps(u_2) + \frac12  g_\eps \int_{\R^2} | u_1|^2 |u_2|^2 \,,
\ee where $E_\eps$ is given by (\ref{enetaeps}) and $g_\eps$ is the intercomponent coupling strength.
The energy is minimized under the mass constraints
 \be  \label{mass}
	\int_{\R^2} |u_j|^2 =   \alpha_j \quad  \text{ with } \quad   \alpha_j>0  \quad  \text{ and } \quad \alpha_1+\alpha_2=1\,.
 \ee

 In \cite{MaAf}, numerical simulations have been
     performed to classify the ground states  according to the values of $\eps$, $g_\eps$  and also the rotational velocity.
     For $\eps$ small and $g_\eps$ large, the numerical evidence is that, for $\alpha_1=\alpha_2=1/2$, the preferred ground state is such that each component is asymptotically located
in a half disk with a local inverted parabola profile. If $\alpha_1 \neq \alpha_2$, they occupy sections in a disk, the area of which is proportional to $\alpha_i$. In particular, when neither $\alpha_i$ is too small, this configuration has less energy than
 a disk vs annulus configuration, which also provides segregation but preserves symmetry. Observation of symmetry breaking
  has also been obtained experimentally very recently \cite{McC}. The breaking of symmetry has been analyzed in \cite{JRL} in a different limit, namely in the case $\eps$ large and $g_\eps $ large.

 Here, we assume strong coupling between components, that is, $g_\eps \to \infty$, and we study the regime

 \be \label{regime}
	g_\eps \eps^2 \to + \infty \quad \text{ and } \quad \eps \to 0\,.
\ee

       A trick introduced in
     \cite{MaAf} is to use a spin formulation
     also called the nonlinear sigma model. In our special setting, since the ground states are  non vanishing real functions, this
     amounts to defining

\be \label{chvariables}
	v := \frac{ \sqrt{|u_1|^2+|u_2|^2}}{\eta_\eps} \, \quad \text{and} \quad \frac\varphi2 := \text{Arg} \left( \frac{|u_1|+i |u_2|}{\sqrt{|u_1|^2+|u_2|^2}}   \right)  \,,
\ee

where $\eta_\eps$ is defined in (\ref{definitionetaeps}). The definition of $\vp$ implies that $|u_1|^2-|u_2|^2=\eta_\eps^2 v^2 \cos{\varphi}$. The mass constraints (\ref{mass}) can be written as

 \be \label{mass2}
	\int_{\R^2} \eta_\eps^2 v^2=\alpha_1+\alpha_2=1 \qquad \text{ and } \qquad \int_{\R^2} \eta_\eps^2 v^2 \cos{\varphi} =\alpha_1-\alpha_2 \,.
\ee We point out that $\cos \vp$ corresponds to the third component of the spin function. Because there is no rotation in the system, the ground states are, up to multiplication by a complex number of modulus one, positive functions. Thus, the second component of the spin is zero and the first one is $\sin \vp$.

 Since the components
     are expected to segregate, the expected behaviour is thus that $v$ tends to 1 except
 on a transition line corresponding to the interface between the two components, while $\varphi$ tends
 to 0 on component 1 and $\pi$ on component 2. This is what we want to analyze rigorously.

 We split the energy into its main contributions and will prove that
\be\label{splitener}
		 \E_\eps(u_1, u_2)  =  E_\eps(\eta_\eps) + F_{\eps}(v) + G_\eps(v,\varphi) \,
 	\ee where $E_\eps$ is given by (\ref{enetaeps}), $\eta_\eps$ is the ground state of $E_\eps$ and
 \ba
	  F_{\eps}(v) &=&   \frac12  \int_{\R^2} \eta_\eps^2 |\nabla v|^2 +  \frac1{2\eps^2} \eta_\eps^4 \{1-v^2\}^2 \,, \label{Feps}\\
	  \label{Geps}  G_\eps(v,\varphi) &=&  \frac18  \int_{\R^2}  \eta_\eps^2 v^2  \, |\nabla \varphi|^2 + \eta_\eps^4 v^4 \, \tilde g_\eps  \{ 1-\cos^2(\varphi)\} \,
\ea
and $ \tilde g_\eps = g_\eps  \Big(1-\frac1{g_\eps \eps^2}\Big) $. Since $\eta_\eps^2$ converges to $\rho$ given by (\ref{rho}) in $\D$, the limits
 of $F_\eps$ and $G_\eps$ can be analyzed as the limits of

 \begin{equation}\label{Fapp} \frac12  \int_{\D} \rho |\nabla v_\eps|^2 +  \frac1{2\eps^2} \rho^2 \{1-v_\eps^2\}^2\end{equation}

 \begin{equation}\label{Gapp} \frac18  \int_{\D}  \rho v_\eps^2  \, |\nabla \varphi_\eps|^2 + \rho^2 v_\eps^4 \, \tilde g_\eps  \{ 1-\cos^2(\varphi_\eps)\}.\end{equation}

 These two energies are of Modica Mortola types with a weight which vanishes on the boundary of $\D$.
  Given the definition of $\varphi_\eps$, there is a domain where $\cos \varphi_\eps$ tends to 1 (asymptotic region of component 1)
   and a domain where $\cos \varphi_\eps$ tends to $-1$ (asymptotic region of component 2), and thus a transition region
    exists between the two domains.
  Two options exist for $v_\eps$: \begin{itemize}\item either $v_\eps$ goes to 1 everywhere, which makes the first energy small
   and the second energy of order $\sqrt {\tilde g_\eps}$, \item or $v_\eps$ goes to zero on the transition line
    where $\cos \varphi$ varies from $+1$ to $-1$: this makes the second energy of lower order and the first energy
     of order $C/\eps$.\end{itemize}
     Because of our hypothesis that $\eps^2\tilde g_\eps$ tends to infinity, it is the second scenario which
     costs less energy. Though $v_\eps$ goes to 1 on each component, it has a transition region of size $\eps$ where it goes sharply to zero. The second energy is of lower order and cannot be seen in the limit. It has just the effect of creating a small region around the interface where $v_\eps$ is small.
     The first energy can be analyzed with techniques coming from \cite{AT} and, once the rescaling in $\eps$ is made, the $\Gamma$-limit comes from the problem on lines:
     \ben
	 I(x)=  \inf \left\{  \frac12  \int_0^\infty  \rho(x) (w')^2 + \frac12 \rho(x)^2 (1-w^2)^2 \,;\, w \in \text{Lip}(\R_+) \,, w(0)=0 \text{ and } w(+\infty)=1 \right\} \,.
\een Using the Euler-Lagrange equation associated with $I$, we shall see that for $x \in \D$, the infimum is attained by the function

\ben
	w^x(t) = \tanh\left( \sqrt{ \frac{\rho(x)}{2} } \,  t\right)  \,,
\een

and we shall have

\be \label{infI}
	 I(x)=\sigma \rho(x)^{\nicefrac32} \hbox{ with }\sigma= \frac1{\sqrt{2}} \int_0^1 \{ 1 - t^2\} \, dt \,.
\ee

This means that $w^x$ is the optimal profile transition at the point $x$, and that $\sigma\rho(x)^{\nicefrac32}$ is the minimum energy needed by $w$, to go from $0$ to $1$ at $x$. In the 1D direction, this provides a weight $2\sigma\rho(x)^{\nicefrac32}$ because as $\eps \to 0$, $v_\eps$ goes from $1$ to $0$ on one side of the interface between the two components, and from $0$ to $1$ on the other side. Therefore, we expect the limit
 to be defined as the integral on the interface where $\vp$ goes from 0 to $\pi$ of the function $2\sigma\rho(x)^{\nicefrac32}$. This requires a precise mathematical definition for this interface.

We define $X$ as the space of functions $\varphi \in BV_{\text{loc}}(\D \,; \{0,\pi\})$ such that

\be \label{masslimit}
	 \int_{\R^2} \rho \cos{\varphi} =\alpha_1-\alpha_2 \,.
\ee

We will prove the $\G$-convergence of  $\eps(\E_\eps(\cdot,\cdot)-E_\eps(\eta_\eps)) $ to $\F$ given in $X$  by

\ben
	\F(\vp) = \frac{2 \sigma}{\pi} \int_{\D} \rho^{\nicefrac32} \, |D\varphi|  \,.
\een

The limiting energy $\F$  measures the length, with a weight of $\rho^{\nicefrac32}$, of the interface between the two phases of $\vp$. Each phase of $\vp$ corresponds to one component of the totally segregated two-component limiting condensate. Notice that when $\F(\vp)$ is finite, $\{ \vp = \pi \}$ has finite perimeter in compact subsets of $\D$, and

\ben
	\F(\varphi) =  2\sigma  \int_{\D\,  \cap \, \partial^*\{\vp = \pi\} }  \rho^{\nicefrac32}  d\H^1 = 2\sigma  \int_{\D\,  \cap \, S\vp }  \rho^{\nicefrac32}  d\H^1 \,.
\een

Here $ \partial^*\{\vp = \pi\} $ stands for the reduced boundary of $\{\vp = \pi\} $ and $S\vp$ is the complement of the Lebesgue points of $\vp$, that is,

\ben
	S\vp = \left\{ x \in \D \,;\, \nexists \, t \in \R \text{ such that } \lim_{r \to 0^+} \frac1{\pi r^2} \int_{B_r(x)}|\vp(y)-t| \, dy = 0 \right\} \,.
\een

We refer to \cite{AFP,EG,Giusti} for the geometric measure theory concepts.  We also refer to \cite{alb} for an introduction to the theory of $\G$-convergence and to the Modica-Mortola theorem by G. Alberti.

We now state our main theorem:
\begin{theo} \label{theo}Let us assume that $V(x) = |x|^2$, and let
\ben
	\mathcal{H} = \left\{ (u_1,u_2) \in H^1(\R^2;\R) \times H^1(\R^2;\R)    \,, \int_{\R^2} V (u_1^2+u_2^2) < \infty \,, (u_1,u_2) \text{ satisfies }   (\ref{mass})   \right\} \,.
\een The functional $\eps(\E_\eps(\cdot,\cdot)-E_\eps(\eta_\eps)) $ $\G$-converges with respect to the $L^1_{loc}(\D) \times  L^1_{loc}(\D)$ distance to $\F(\vp)$, in the following sense:
	
	\texttt{(Compactness)} for every sequence $\{(u_{1,\eps}, u_{2,\eps})\}_{\eps>0}$  of minimizers of $\E_\eps$ in $\mathcal{H}$ such that
	
	\be \label{upgcvtheo}
		\sup_{\eps >0} \eps \left( \E_\eps (u_{1,\eps},u_{2,\eps} ) - E_\eps(\eta_\eps) \right) < +\infty \,,
	\ee
there exists  $\varphi \in X$ and a (not relabeled) subsequence such that

\ba
	(u_{1,\eps},u_{2,\eps} ) \to \sqrt{\rho}   \left( \mathbf{1}_{\{ \vp = 0\}} ,  \mathbf{1}_{\{ \vp = \pi \}}  \right) \quad &\text{ in }& \quad L^1_{loc}(\D) \times  L^1_{loc}(\D)  \,; \label{compactnessgcvtheo}
\ea
and \texttt{(Lower bound inequality)}

\be  \label{lowerboundgcvtheo}
	\liminf_{\eps \to 0}  \eps \left( \E_\eps (u_{1,\eps},u_{2,\eps} ) - E_\eps(\eta_\eps) \right) \geq \F(\vp) \,.
\ee	
	\texttt{(Upper bound inequality)}   For every $\vp \in X$, there exists a sequence $\{(u_{1,\eps}, u_{2,\eps})\}_{\eps>0} \subset \mathcal{H}$, converging as  $\eps \to 0$ to $ \sqrt \rho\left( \mathbf{1}_{\{ \vp = 0\}} ,  \mathbf{1}_{\{ \vp = \pi \}}  \right) $ in $L^1_{\text{loc}}(\D) \times L^1_{\text{loc}}(\D)$, such that

\be  \label{upperboundgcvtheo}
	\limsup_{\eps \to 0} \eps \left( \E_\eps (u_{1,\eps},u_{2,\eps} ) - E_\eps(\eta_\eps) \right)  \leq \F(\vp)  \,.
\ee
\end{theo}

We point out that we only prove the $\G$-convergence at the level of minimizers of $\E_\eps$. Indeed, minimizers of the functional have the property that they are positive functions which do not vanish. Therefore, this property allows the definition of $(v,\vp)$ through (\ref{chvariables}). As usual, the $\G$-convergence theorem implies the convergence of the energy of the ground states:

\begin{coro} \label{coro}

If $\{(u_{1,\eps}, u_{2,\eps})\}_{\eps>0}$ is a sequence of minimizer of $\E_\eps$ in $\mathcal{H}$, then

\be \label{limitenergy}
	\lim_{\eps \to 0} \, \eps \left(  \E_\eps(u_{1,\eps}, u_{2,\eps}) -  E_\eps(\eta_\eps) \, \right) = \inf_X \F \,.
\ee
\end{coro}

A study of the ground states of $\F$ allows us to prove symmetry breaking when neither $\alpha_i$ is too small:
\begin{coro} \label{corobreak}

There exists $\delta_0$ of order 0.15, such that if $\alpha_1 \in [\delta_0,1-\delta_0]$, then for $\varepsilon$ sufficiently small,  the  minimizers $(u_{1,\eps}, u_{2,\eps})$ of $\E_\eps$ in $\mathcal{H}$ are not radial.
\end{coro}

\begin{rema}Our main theorem remains true when $V$ is any trapping potential for which we have good estimates for the ground state $\eta_\eps$, namely the estimates in Proposition \ref{estimatesetaeps}.\end{rema}

\subsection{Links with related problems}
 The segregation behaviour in two component condensates has been widely studied: regularity of the wave function \cite{ctv3,NoTaTeVe,WeWe1}, regularity of the interface \cite{CaffLin2}, asymptotic behaviour near the interface \cite{BeLinWeiZhao,BeTer}.
 The main difference with these references is that, on the one hand, we use mainly the energy instead of the equation and, on the other hand,  we do not switch off the trapping potential by blowing up the problem near the interface or by
  considering a bounded domain with no trapping. Indeed, we consider the limit where $\eps$ goes
  to zero at the same time as $g_\eps \eps^2$ going to infinity, so that it is the trapping potential which provides the leading order  behaviour of the wave function through the inverted parabola profile $\rho$. In all the previous quoted references,
   $\eps$ is set to 1, so that in the limit $g_\eps$ large,
  the trapping potential is not present, and the limiting profile is 1. We deal with the trapping potential by a proper division of the limiting wave function which allows to express nicely the energy using a trick introduced by \cite{LM}. Nevertheless,
   our proofs which rely on energy considerations also provide information for the case $\rho=1$.

In [31], the authors fix a point $x_\infty$ on the interface $\partial A$,
and consider a sequence $x_\eps$ tending to $x_\infty$ such that
 $u_{1,\eps}(x_\eps)=u_{2,\eps}(x_\eps)=m_\eps$. An open question in [31]  is to prove
 in 2D that $g_\eps m_\eps^4$ stays bounded. This may be obtained
 with our technique since  in our case $m_\eps$
 is probably related to the minimum of $v_\eps$. We detail this remark in Section 5.3.

\subsection{Main ideas in the proof}
Let us now give more details on the proof.

The proof consists of upper and lower bounds, that we construct for the functional $\mathcal F_\eps (v_\eps, \vp_\eps)= \eps \left(  \E_\eps(u_{1,\eps}, u_{2,\eps}) -  E_\eps(\eta_\eps) \, \right)$.

For the upper bound, we choose the set $A$ where asymptotically $u_2$ will be $\rho$.
 In a first step, we assume that  $\vp = \pi \mathbf{1}_A$, where $A$ is an open bounded subset of $\R^2$ with smooth boundary such that $\H^1(\pA \cap \pD)=0$.
  The test function $\vp_\eps$ is matched between 0 in a subdomain of $\D \setminus \bar A$ to $\pi$ in a subdomain of $A$, using a transition region of size $\eps t_\eps$. In order to approximate the optimal 1 dimensional profile that solves $I(y)$, we define
  \ben
	w_{\eps,T}  = \left\{ \begin{array}{cll} m_\eps &   \text{ in } &  (0,  t_\eps) \\
  \tanh &   \text{ in } &  (t_\eps,  T) \\
		h &   \text{ in } &   (T,T+\nicefrac1T)    \\ 		
		1 &   \text{ in }  & (T+\nicefrac1T,+\infty ) , \quad
\end{array} \right.
\een  where $t_\eps=\tanh m_\eps$ and  $h$ is a polynomial which matches smoothly $\tanh$ to 1. Then we define \ben
	w_{\eps,T}^y(t) = w_{\eps,T}\left( \sqrt{\frac {\rho(y)}2} \, t\right) \,,
\een for $t= d(x)/\eps <CT$, and $d(x)$ is the distance to the boundary. In order to construct $v_\eps$, we need a partition
 of unity for $\partial A$, where we match the functions $w_{\eps,T}^{y_i}$, as $y_i$ varies along this partition. For this
  $v_\eps$, we can estimate $F_\eps$ with techniques similar to those of Modica Mortola \cite{Mod}, and to the adaptation
   of these techniques to problems with weight by Bouchitté \cite{BOU}. Because $\rho$ vanishes, we cannot use directly the results of Bouchitté and we need precise estimates on the behaviour of $\eta_\eps$ near the boundary.
   Since $w_{\eps,T}$ is the optimal profile for the 1D version of (\ref{Fapp}),  there is a transition from 1 to 0 and a transition from 0 to 1 and we find an upper bound which is $2\int_{\partial A}I(y)\ dy$.
  Then we prove that for this test function, $G_\eps (v_\eps,\vp_\eps)$ is lower order: indeed, the transition layer for $\vp$ is is of order $\eps t_\eps$, so much smaller than the one of $v_\eps$. Hence in $G_\eps$, $v_\eps$ can be approximated by $m_\eps$. We choose $m_\eps^4=\eps^2 g_\eps$, which tends to 0, and makes $G_\eps$ of lower order.

This provides the upper bound for an open bounded subset $A$ with smooth boundary such that $\H^1(\pA \cap \pD)=0$.
  We show in the appendix that for any $\vp\in X$, $\{\vp=\pi \}$ can be approximated by  sets $A$ which are open bounded subsets of $\R^2$ with smooth boundary such that $\H^1(\pA \cap \pD)=0$ and that the mass constraints can be satisfied for the approximating $u_{1,\eps}, u_{2,\eps}$.

  The difficulty in the lower bound is to prove that $v_\eps$ goes to zero on a line and that it provides a positive lower bound. Indeed, the usual Modica-Mortola bound would imply that $v_\eps$ goes to 1 almost everywhere and the lower bound is 0. We have to use $G_\eps$ and the upper bound to prove that
  $v_\eps$ has a transition to 0 and that $\cos ^2\vp_\eps$ tends to 1. Hence, because of the mass constraint, we get two regions where asymptotically $\vp_\eps $ is 0 and $\pi$. To analyze the behaviour of $v_\eps$, we use  the slicing method introduced in \cite{AT} (see also \cite{braides}). This consists in looking at the transition for $v_\eps$ in one dimensional slices and get the 1D energy estimate.  The use of the energy $G_\eps$ is only to prove that $v_\eps$ goes to zero.   We first prove the lower bound for $\eps F_\eps$ in 1D using the coarea formula, and then in 2D using
   the slicing method. We get that  $\eps F_\eps(v_\eps,\vp_\eps,  E)$ converges to a measure $\mu(E)$ supported in $S_\vp$ of density $\rho^{3/2}$ with respect to the $\H^1$ measure. The last part of the proof of the lower bound is inspired by ideas in \cite{ABS}.

We end  with a variant of the coarea formula that can be found in \cite{DM} Lemma 2.2, and in \cite{BOU} Proposition 2.

\begin{prop} \label{coareaprop}
	Let $\Omega$ be an open bounded subset of $\R^N$, and $\Psi(x,s,p)$ a Borel function of $\Omega \times \R \times \R^N$, which is sublinear in $p$. Let $u$ be a Lipschitz continuous function on $\Omega$ and denote, for every $t>0$, $S_t = \{ x \in \Omega \,;\,u(x) <t \}$. Then, for almost every $t \in \R$ , $\mathbf{1}_{S_t}$ belongs to $BV(\Omega)$ and we have
	
\be \label{coarea}
	\int_\Omega \Psi(x,u,Du) \, dx  = \int_{-\infty}^{\infty} dt \int_{\Omega}  \Psi(x,t,D\mathbf{1}_{S_t}) \,.
\ee	
	
\end{prop}

The paper is organized as follows: in Section 2, we present the properties of $\eta_\eps$. Then in Section 3, we prove the
 decoupling of energy (\ref{splitener}) and how to go from the $(u_1,u_2)$ formulation to $(v,\vp)$. Section 4 is devoted to the upper bound, and Section 5 to the lower bound. Finally, in Section 6, we prove our main theorem.

%%%%%%%%%%%%%%%%%%%%%%%%%%%%%%%%%%%%%%%%%%%%%%%%%%%%%%%%%%%%%
\subsection{To go further}

\subsubsection{Analysis of the limiting problem}
A natural question is to analyze the limiting problem, that is the ground state of ${\cal F}$ under the constraint (\ref{masslimit}).
 If we define $A$ to be the set where $\cos \varphi =1$. Then $\int_A \rho=\a_1$ and $\int_{\D\setminus A} \rho =\a_2$ with $\a_1+\a_2=1$.

 If $\rho=1$, then the problem of minimizing ${\cal F}$
 amounts to minimizing $|\partial A|$ under the constraints $|A|=\a_1 $ and $|\D\setminus A|=\a_2=1-\a_1 $. The Euler-Lagrange equation
  of the minimization problem yields that the curvature is either 0 or constant, hence $A$ is either a disk, an annulus or a disk sector. The equivalent problem with a weight $\rho$  is open.

   If we assume that the solution is either two disks sectors or a disk and an annulus, we can compute
   explicitly the energy ${\cal F}$ and find  that if $\a_1=\a_2$, then the optimal configuration is two half disks, while
   if $\a_1$ is much less then $\a_2$, then the ground state is a disk and an annulus (see Section 6.4). Indeed, the energy of two disk sectors is $3\sigma/2$, while the energy of a disk and annulus is $8\sigma (1-\alpha_1)^{\nicefrac34} (1-\sqrt{1-\alpha_1})^{\nicefrac12} $ if $\alpha_1$ corresponds to the mass of the inside disk.
 If $\alpha_1$ or $\alpha_2=1- \alpha_1$ is to small, then the disk and annulus becomes the preferred configuration.
   In the case $\alpha_1=\alpha_2=1/2$, it follows from our theorem that symmetry breaking occurs since at the limit, the disk plus annulus configuration does not minimize the energy. These two cases are well illustrated in the experimental observations of \cite{McC}, figure 4.

We insist on the point that a rigorous analysis of the ground states of $\mathcal F$ in $X$ is an interesting open question.

\subsubsection{Convergence for $u_{1,\eps}$, $u_{2,\eps}$}
 The convergence that we have for $(u_{1,\eps}$, $u_{2,\eps})$ to $\sqrt \rho( \mathbf{1}_{\{ \vp = 0\}} ,  \mathbf{1}_{\{ \vp = \pi \}}  )$ is very weak. Nevertheless, we expect that on compact subsets of $\mathbf{1}_{\{ \vp = \pi \}}$ or
  $\mathbf{1}_{\{ \vp = 0\}}$, the convergence can be improved. For instance, it would be natural to have similar convergence as that of $\eta_\eps$ to $\sqrt \rho$ (that is $C^1_{loc}$) on these domains.

  \subsubsection{Case $g_\eps \eps^2$ of order 1}
An interesting open question is to deal with the case when $g_\eps \eps^2$ tends to a positive finite constant $c_0^2$.
 In this case, $F_\eps$ and $G_\eps$ become of the same order and we expect that $m=\liminf_{\eps\to 0} v_\eps$ is  a positive constant (on the interface where $\vp$ varies), instead of being 0. We believe that our techniques
  still provide an upper bound for the problem. We expect the $\Gamma$ limit to be
  $$\left ( {2 \sigma_m}+c_0 \frac \pi 4 m^3\right ) \frac 1 \pi \int_{\D} \rho^{\nicefrac32} \, |D\varphi|  \,.$$
   where $\sigma_m=\frac 1 {\sqrt 2} \int_m^1 (1-t^2)\ dt$.

 \subsubsection{Case of different scattering lengths}

 In this paper, we consider that the scattering lengths are the same for both components,
 that is, in (1.4) it is the same energy $E_\eps$  for both components.
 When the two components result experimentally from different atoms, the
two scattering lengths are very close but not equal.
 This leads to an energy $E_\eps$ depending on the component, namely

    $$E_{\eps, i}(\eta) = \frac12 \int_{\R^2} |\nabla \eta|^2 + \frac1{\eps^2} |x|^2 |\eta|^2  +   \frac{g_i}{2\eps^2} |\eta|^4   \,,$$

    where $g_i$ is related to the scattering lentght of component $i$. If $g_1\neq g_2$, then  the leading order Thomas Fermi approximation is no longer
    the same for each component, namely it is
    $$g_i \rho_i=\lambda_i^2-|x|^2 \hbox{ in } B_i=B(0,\lambda_i).$$
     The limiting problem becomes: find a partition of $B_1 \cup B_2$ into three sets $A_1$, $A_2$ and $N$,
     such that $u_{i,\eps}^2 \to \rho_i \mathbf{1}_{A_i} $, $\int_{A_i} \rho_i  =\alpha_i$ and it minimizes
     \begin{equation}\label{rhoi}\int_{A_1} |x|^2 \rho_1+\frac{g_1}2 \rho_1^2 + \int_{A_2} |x|^2 \rho_2+\frac{g_2}2 \rho_2^2.\end{equation}
     This problem is open and is probably related to the problem of finding a partition of the disk into
     two subdomains which minimize the sum of the first eigenvalues of the Dirichlet laplacian.

     Of course, in our case, since we have $B_1=B_2$, $\rho_1=\rho_2$ and $N=\emptyset$, (\ref{rhoi})  does not provide
     any information at leading order. This is why we have to go to the next order which yields the perimeter
      minimization problem.

%%%%%%%%%%%%%%%%
%%%%%%%%%%%%%%%%

%%%%%%%%%%%%%%%%
%%%%%%%%%%%%%%%%
\section{Estimates for $\eta_\eps$} \label{estimates}

Let $\eta_\eps$ be the ground state defined by (\ref{definitionetaeps}). The ground state is a non vanishing radially symmetric function. It is unique up to multiplication by a constant of modulus one, and satisfies  the Gross-Pitaevskii equation

\be \label{eqetaeps}
	 - \Delta  \eta_\eps + \frac1{\eps^2} |x|^2 \eta_\eps + \frac1{\eps^2}  |\eta_\eps|^2\eta_\eps = \frac{\lambda_\eps}{\eps^2} \, \eta_\eps \,.
\ee

The term $\eps^{-2}\,\lambda_\eps$ is the Lagrange multiplier associated with the mass constraint, and the pair $(\eta_\eps, \lambda_\eps)$ is unique among positive solutions of (\ref{eqetaeps}). As $\eps$ tends to 0, $\eta_\eps$ tends to $\sqrt \rho$ given by (\ref{rho}). Throughout the paper, we will need precise estimates for this convergence. The following proposition, based on previous results in \cite{AJR,Ga,GaPe,IM,KaSou}, sums up the properties of $\eta_\eps$. We point out that it follows from \cite{Ga,GaPe,KaSou} that an approximation of $\eta_\eps$ by $\sqrt \rho$ holds as close to the boundary of $\D$ as needed and is given by (\ref{estetanearboundary}).  We also include an estimate of $\rho$ in terms of the distance to the bulk that will be used in the proofs.

\begin{prop} \label{estimatesetaeps} There are constants $c,C>0$, $\alpha \in (\nicefrac12,\nicefrac35)$ and $\gamma \in (\nicefrac12,\nicefrac34)$, such that for $\eps$ sufficiently small, $\rho,\lambda$ being given by (\ref{rho}),

\ba
%	E_\eps(\eta_\eps) &\leq& C \, \eps^{-2}  \label{boundEeps} \, \\
	E_\eps(\eta_\eps) &\leq& C/\eps^2  \,\label{boundE1eps}, \\
	| \lambda_\eps - \lambda | &\leq& C \, \eps \, |\ln \eps|^{\nicefrac12} \, \label{lambdaeps}, \\
	\| \eta_\eps - \sqrt{\rho}\|_{C^1(K)} &\leq& C_K \, \eps^2 \, |\ln \eps| \hspace{1.425cm} \text{ for  }  K \subset\subset \D  \label{estetaK} \,, \\
	| \eta_\eps(x)  - \sqrt{\rho}(x) | &\leq& C \, \eps^{\gamma} \hspace{2.65cm} \text{ for } x \in B(0\,,\lambda- c\, \eps^\alpha)  \label{estetanearboundary} \,, \\
	\eta_\eps(x) &\leq& C \, \eps^{\nicefrac16} \, e^{c\,\eps^{-\nicefrac 13}(\lambda- |x|)} \hspace{.3cm}  \text{ for } x \in \R^2 \backslash \D \label{estetaoutbulk} \,, \\
	\partial_r \eta_\eps(|x|) &\le &0   \hspace{3.2cm}  \text{ for } x \in \R^2, \label{estetarad} \\
	\frac{\rho(x)}{\lambda \, dist(x, \pD) } &\in& [1,2)  \hspace{2.55cm} \text{ for } x \in \D  \label{estrhopD} \,.
\ea 	
\end{prop}

\textbf{Proof:} for the proof of (\ref{boundE1eps}), one can rewrite the energy as

\be \label{EE1}
	E_\eps(\eta) = E^1_\eps(\eta) + \frac1{2\eps^2} \left( \lambda^2 -  \frac12 \int_{\D} \rho^2 \right)
\ee

where

\ben
	E^1_\eps(\eta) =  \frac12 \int_{\R^2} |\nabla \eta|^2 + \frac1{2\eps^2} \left(|\eta|^2 - \rho(x) \right)^2  +   \frac1{\eps^2} (\lambda^2 - |x|^2)_- \,  |\eta|^2   \,,
\een
 and $(\lambda^2 - |x|^2)_-$ is the negative part of $(\lambda^2 - |x|^2)$. In Theorem 2.1 of \cite{AJR},
 it is proved that $E^1_\eps(\eta)\leq C |\ln \eps |$. Then (\ref{boundE1eps}) follows from (\ref{EE1}) and the fact that
 $\int_{\D} \rho^2 =2\lambda^2/3$.

 Estimate (\ref{estetaK}) is proved in Proposition 2.2 of \cite{IM}.
Estimates (\ref{lambdaeps}) and  (\ref{estetaoutbulk}) are proved in Theorem 2.1 of \cite{AJR}. Estimate (\ref{estetarad}) is also proved in Theorem 2.1 of \cite{AJR}, but only in a neighborhood of $\partial\D$. But the proof, however, works in the case $V(x)=|x|^2$ and the estimate holds in all $\R^2$.\\

We now prove (\ref{estetanearboundary}). For $\lambda>0$, we define $\tilde \eta_{\eps,\lambda}$ as the unique radially symmetric, positive solution  of the equation

\be \label{eqnmcl}
		 - \eps^2 \Delta  \eta +  (\lambda^2-|x|^2) \eta +  \eta^3  = 0 \,.
\ee

The function $\tilde \eta_{\eps,\lambda}$ corresponds to a ground state of a BEC without mass constraint. In \cite{Ga,GaPe,IM}, the behavior of $\tilde \eta_{\eps,\lambda}$ is studied. Using the results in  Proposition 1.2, Remark 1.3 and Proposition 1.4 in \cite{Ga}, we obtain

\ben
	\tilde \eta_{\eps,1}(x) = \eps^{\nicefrac13} \nu_0\left( \frac{1-|x|^2}{\eps^{\nicefrac23}}\right) + \O(\eps) \,,
\een
	
where

\ben
	\nu_0(y) = y^{\nicefrac12} - \frac12 y^{-\nicefrac52}+  \O_{y \to +\infty}(y^{-\nicefrac{11}2}) \,, \quad y \in (-\infty, \eps^{-\nicefrac23}] \,.
\een

Hence, for $x \in B(0,1)$ we obtain

\ben
	|\,\tilde \eta_{\eps,1}(x) - \sqrt{1-|x|^2} \, | \leq C \,   \left(  \eps^2 (1-|x|^2)^{-\nicefrac52} + \eps^4 (1-|x|^2)^{-\nicefrac{11}2} + \eps \right) \,.
\een

In particular, if $x \in  B(0\,,\lambda- \eps^\alpha) $ with $\alpha \in (\nicefrac12,\nicefrac35)$, we get

\be \label{estetanmc}
	|\,\tilde \eta_{\eps,1}(x) - \sqrt{1-|x|^2} \, | \leq C \,   \left(  \eps^2 \, \eps^{-\nicefrac{5\alpha}2} + \eps^4 \, \eps^{-\nicefrac{11\alpha}2} + \eps \right) = \O(  \eps^\gamma ) \,
\ee

with $\gamma \in  (\nicefrac12,\nicefrac34)$. We will use (\ref{estetanmc}) to prove (\ref{estetanearboundary}). First, a straight computation shows that defining $\eps_{\lambda} = \lambda^{-2} \eps$, $\tilde \eta_{\eps_{\lambda},\lambda}$ solves equation (\ref{eqnmcl}) with $\lambda=1$. Hence, considering (\ref{estetanmc}), a change of variables gives

\be \label{estetanmc2}
	|\,\tilde \eta_{\eps_{\lambda},\lambda}(x) - \sqrt\rho(x) \, |  = \O(  \eps^\gamma ) \,,
\ee

for  $x \in B(0\,,\lambda- (\lambda^{-2} \eps)^\alpha)$. In Proposition 2.2 and Theorem 2.2 in \cite{IM}, it is proved that

\be \label{derivative}
	\| \nabla \eta_\eps \|_{L^\infty(\R^2)} = \O(\eps^{-1}) \,;
\ee

and that

\be \label{cheta}
	\eta_{\eps, \lambda}(x) = \ell_{\eps,\lambda}^{\nicefrac12} \,\, \tilde \eta_{\tilde \eps, \lambda}(\ell_{\eps,\lambda}^{-1} \,\, x) \,,
\ee

where

\ben
	\ell_{\eps,\lambda} = \left( 1 +  \frac{\eps \lambda_\eps }{\lambda} \right)  \quad \text{ and }  \quad \tilde \eps = \ell_{\eps,\lambda}^{-1} \eps \,.
\een

It follows from (\ref{lambdaeps}) that

\ben
	\ell_{\eps,\lambda} =  1+\O(\eps^2 \,|\ln \eps|^{\nicefrac12})   \quad \text{ and } \quad \tilde \eps =  \eps +\O(\eps^2 \,|\ln \eps|^{\nicefrac12}) \,.
\een

Hence, using  (\ref{derivative}) and (\ref{cheta}), we obtain

\ben
	\eta_{\eps, \lambda}(x) = \tilde \eta_{ \tilde \eps, \lambda}(x)  + \O(\eps \,|\ln \eps|^{\nicefrac12}) \,.
\een

Putting this last estimate in (\ref{estetanmc2}), and using that $\gamma \in (\nicefrac12,\nicefrac34)$, we obtain

\ben
	|\, \eta_{\eps_{\lambda}(x),\lambda}- \sqrt\rho(x)\, |  = \O(  \eps^\gamma) \,,
\een

for $x \in B(0\,,\lambda- c \, \eps^\alpha)$ with $c>0$. We derive (\ref{estetanearboundary}) by changing $\eps_\lambda$ by $\eps$ in the previous estimate. Finally, writing

\ben
		\frac{\rho(x)}{\lambda \, dist(x, \pD) } =  \frac{(\lambda + |x|)}{\lambda}
\een

we get  (\ref{estrhopD}) for $|x| < \lambda$. \\ \qed

%%%%%%%%%%%%%%%%
%%%%%%%%%%%%%%%%
\section{Rewriting the energy} \label{energyrewriting}

In this section, we prove equality (\ref{splitener}), that is, the reformulation of the Gross-Pitaevskii energy of a two component condensate in (\ref{enertwobec}), as the weighted  Cahn-Hilliard energy for the pair $(v,\vp)$ defined by  (\ref{chvariables}), plus the energy of the ground state $\eta_\eps$ of a one component condensate. We start by giving the properties of the minimizers of $\E_\eps$ and the properties of the corresponding pairs $(v_\eps,\vp_\eps)$ defined by (\ref{chvariables}).

%%%%%%%%
\begin{prop} \label{propu1u2vphi} \textbf{(i)}  Let $\{(u_{1,\eps},u_{2,\eps})\}_{\eps >0}$ be a sequence of minimizing pairs of  $\E_\eps$ in $\mathcal{H}$ satisfying (\ref{upgcvtheo}). Then, each component is a non vanishing smooth function, and there is $C>0$ such that

 \be \label{booundlinftyu1u2}
 	\|u_{1,\eps}\|_{L^\infty(\R^2)} \,, \|u_{2,\eps}\|_{L^\infty(\R^2)} < C
 \ee

 for every $\eps>0$. Moreover, the pairs $(v_\eps,\vp_\eps)$ are well defined by (\ref{chvariables}), verify the mass constraints (\ref{mass2}) and we have

 \be \label{vphilocLipschitz}
 	(v_\eps,\vp_\eps) \in Lip_{loc}(\R^2 \,;\, (0,+\infty) \times [0,\pi])
\ee
and
\be \label{unifupperboundveps}
	\sup_{\eps>0} \|v_\eps\|_{L^{\infty}(K)} < C_K \qquad \text{ for every } K \subset \subset \D \,.
\ee

\textbf{(ii)} Conversely, let $(v,\vp) \in Lip(\R^2 \,;\, (0,+\infty) \times [0,\pi])$ satisfying (\ref{mass2}) such that $v\,, \nabla v \,, \nabla \vp \in L^{\infty}(\R^2)$. Then, defining

\be
	u_{1} =  \, \eta_\eps \,  v \,  \cos{(\vp/2)}    \quad \text{ and } \, \quad u_{2} =  \,  \eta_\eps v \,   \sin{(\vp/2)}     \, \label{invert} \,,
\ee
we have $(u_{1},u_{2}) \in \mathcal{H}$ and $|u_{1}|^2+|u_{2}|^2>0$.
\end{prop}

\textbf{Proof:}  \textbf{(i)} Let $(u_{1,\eps}, u_{2,\eps})$ be a minimizer of $\E_\eps$ in $\mathcal{H}$.
 Since $\E_\eps ( |u_{1,\eps}|, |u_{2,\eps}|)\leq \E_\eps (u_{1,\eps}, u_{2,\eps})$, the  pair of the absolute values satisfies the system

\ba
	- \Delta u_{1,\eps}  + \big( \eps^{-2} V +  \eps^{-2}  u_{1,\eps} ^2 + g_\eps u_{2,\eps}^2 \big) \, u_{1,\eps} &=& \lambda_{1,\eps } \, u_{1,\eps}  \label{eqELu1} \\
	- \Delta u_{2,\eps}  + \big( \eps^{-2} V +  \eps^{-2}  u_{2,\eps} ^2 + g_\eps u_{1,\eps}^2 \big)  \, u_{2,\eps} &=& \lambda_{2,\eps } \, u_{2,\eps}  \,, \label{eqELu2}
\ea

where $\lambda_{1,\eps }$ and $\lambda_{2,\eps }$ are the Lagrange multipliers associated with (\ref{mass}). The strong maximum principle yields that  $|u_{1,\eps}| $ and $|u_{2,\eps}| $ are positive functions. Using standard elliptic regularity, we deduce further that $u_{1,\eps} $ and $u_{2,\eps} $ are non vanishing smooth functions. We use an argument in \cite{IM} to prove that $u_{1,\eps}$ and $u_{2,\eps}$ are uniformly bounded in $\R^2$. Let us define $w = \eps^{-1} |u_{1,\eps}| -  \lambda_\eps^{\nicefrac12}$. We have $w \in L_{loc}^3(\R^2)$ and $\Delta w \in L_{loc}^1(\R^2)$. Kato's inequality and equation (\ref{eqELu1}) give

\ben
	\Delta (w^+) \geq \text{sgn}^+(w) \,\Delta w \geq\eps^{-3} \, \text{sgn}^+(w) \,\eps w \,(\eps w + \eps \lambda_\eps^{\nicefrac12})\,(\eps w + 2\eps \lambda_\eps^{\nicefrac12}) \geq (w^+)^3 \,.
\een

Hence, $-\Delta (w^+) + (w^+)^3\leq 0$ weakly in $\R^2$ and Lemma 2 in \cite{brezis1} yield $w^+\leq0$. We obtain $ |u_{1,\eps}| \leq \eps \lambda_\eps^{\nicefrac12}$. Multiplying equation (\ref{eqELu1}) by $u_{1,\eps}$ and then integrating we find $ \lambda_\eps^{\nicefrac12} \leq 2 \, \sqrt{\E_\eps ( |u_{1,\eps}|,| u_{2,\eps}|)  }$. Since $( u_{1,\eps}, u_{2,\eps})$ verifies (\ref{upgcvtheo}), from estimate (\ref{boundE1eps}), we derive

\ben
	0< |u_{1,\eps}| \leq \eps \sqrt{\Big( \E_\eps( u_{1,\eps}, u_{2,\eps}) - E_\eps(\eta_\eps) \Big)+E_\eps(\eta_\eps) } < C \,.
\een

We similarly prove that $0<|u_{2,\eps}| < C$, so (\ref{booundlinftyu1u2}) is proved. Since $\eta_\eps>0$, $u_{1,\eps}$ and $u_{2,\eps}$ do not vanish in $\R^2$, the pairs $(v_\eps,\vp_\eps)$ are well defined by (\ref{chvariables}) and $v_\eps >0$. Since $u_{1,\eps}$ and $u_{2,\eps}$ are smooth, $v_\eps$ and $\vp_\eps$ are locally Lipschitz functions so (\ref{vphilocLipschitz}) holds. The definition of $v_\eps$ and (\ref{mass}) give

\ben
	\int_{\R^2} \eta_\eps^2 v_\eps^2 = \alpha_1 +\alpha_2 \,.
\een

From the definition of $\vp_\eps$, we infer that

\be \label{relationvpu12}
	\cos( \vp_\eps ) = \frac{|u_{1,\eps}|^2-|u_{2,\eps}|^2 }{|u_{1,\eps}|^2+|u_{2,\eps}|^2 } \,,
\ee

which, together with (\ref{mass}), yields

\ben
	\int_{\R^2} \eta_\eps^2 v_\eps^2 \cos \vp_\eps = \alpha_1 -\alpha_2 \,.
\een

Hence, $(v_\eps,\vp_\eps)$ satisfies (\ref{mass2}). Finally, the estimate (\ref{estetaK}) gives $\eta_\eps \geq c_K>0$ in $K \subset \subset \D$, so (\ref{booundlinftyu1u2}) yields (\ref{unifupperboundveps}).  \\

\textbf{(ii)} Consider $(v,\vp)$ as in the statement and define $(u_1,u_2)$ by (\ref{invert}). Since $(v,\vp)$ verifies (\ref{mass2}), relation (\ref{invert}) gives

\ben
	\int_{\R^2} |u_{1}|^2 + |u_{2}|^2 = \int_{\R^2} \eta_\eps^2 v^2  = \alpha_1 + \alpha_2
\een	

and

\ben
	\int_{\R^2} |u_{1}|^2 - |u_{2}|^2 =  \int_{\R^2} \eta_\eps^2 v^2 \cos^2 \vp = \alpha_1 - \alpha_2 \,.
\een	

Thus, $(u_{1},u_{2})$ verifies (\ref{mass}). We have $|u_{1}|^2+|u_{2}|^2>0$. Indeed, if it was not the case, since $v>0$ then $\vp$ should take simultaneously the values $0$ and $\pi$. Since $v \in L^\infty(\R^2)$, bounds (\ref{estetanearboundary}) and (\ref{estetaoutbulk}) on $\eta_\eps$ give

\ben
	\int_{\R^2} V \left( |u_{1}|^2 + |u_{2}|^2 \right) \leq C \int_{\R^2} V \eta_\eps^2 < + \infty \,.
\een	

We compute

\ben
	|\nabla u_{1}|^2 \leq   C \left( v^2 |\nabla \eta_\eps|^2 + \eta_\eps^2 |\nabla \eta_\eps|^2 + v^2  \eta_\eps^2   |\nabla \vp|^2  \right)  \,.
\een

The right hand side of the inequality is integrable in $\R^2$ because $v\,, \nabla v \,,  \nabla \vp \in L^\infty(\R^2)$ and $\eta_\eps \in H^1(\R^2) \cap L^\infty(\R^2)$. Thus, $u_{1,\eps} \in H^1(\R^2)$. We prove similarly that $u_{2,\eps} \in H^1(\R^2)$. We have proved that $(u_{1,\eps},u_{2,\eps}) \in \mathcal{H}$. \\ \qed \\

We now prove the rewriting of the energy.

\begin{prop} \label{energysplittingprop} Let $(u_1,u_2) \in \mathcal{H}$ satisfying $|u_1|^2 + |u_2|^2 >0$. Defining $(v,\vp)$  by  (\ref{chvariables}) we have
	\ben
		 \E_{\eps}(u_1, u_2)  =  E_\eps(\eta_\eps) + F_{\eps}(v) + G_{\eps}(v,\varphi) \,,
 	\een
where $E_\eps$, $F_{\eps}$ and  $G_{\eps}$ are given respectively by (\ref{enetaeps}), (\ref{Feps}) and (\ref{Geps}).

\end{prop}

\textbf{Proof:}  since $|u_1|^2 + |u_2|^2 >0$, the pair $(v,\vp)$ is well defined. The definitions of $v$ and $\varphi$  yield

\be  \label{a1}
	|u_1| =  \eta_\eps \, v \, \cos{(\vp/2)} \qquad \text{ and } \qquad |u_2| =  \eta_\eps \, v \, \sin{(\vp/2)}\,,
\ee

which give

\ba
	  |u_1|^2+|u_2|^2  &=&  \eta_\eps^2 v^2 \nonumber \\
	 |u_1|^2|u_2|^2 &=&\frac 14   \eta_\eps^4 v^4 \{ 1 - \cos^2 \varphi \} \label{d} \\
	 |u_1|^4+|u_2|^4 &=& \frac 12   \eta_\eps^4 v^4 \{ 1 + \cos^2 \varphi \} \nonumber \,.
\ea

Since $u_1$ and $u_2$ are real and do not change sign, we have $|\nabla u_1|^2=|\nabla |u_1||^2$ and $|\nabla u_2|^2=|\nabla |u_2||^2$. The relations in (\ref{a1}) give then

\be \label{derivuu}
		|\nabla u_1|^2 + |\nabla u_2|^2 =   |\nabla (v\eta_\eps)|^2  + \frac 1{4}  (v\eta_\eps)^2  |\nabla  \varphi |^2 \,.
\ee

Replacing  (\ref{d}) and (\ref{derivuu}) in $ \E_{\eps}(u_1, u_2) $ we get

\ba \label{Eetaepsv}
	 \E_{\eps}(u_1, u_2)  &=&   \frac12 \int  |\nabla (v\eta_\eps)|^2 + \frac1{\eps^2} V \eta_\eps^2 v^2 \\
	 &+& \frac12 \int   \frac 1{4}  v^2\eta_\eps^2  |\nabla  \varphi |^2+  \frac1{4\eps^2}  \eta_\eps^4 v^4 \{ 1 + \cos^2 \varphi \} + \frac 1{4}   g_\eps  \, \eta_\eps^4 v^4 \{ 1 - \cos^2 \varphi \} \nonumber \,.
 \ea

The previous formulation of the energy is the one given by the spin formulation (see the introduction and \cite{MaAf}). We now show how the phase transition model is obtained. Performing an integration by parts, using (\ref{eqetaeps}) and the first mass constraint in (\ref{mass2}), we obtain

\ba
	 \int  |\nabla (v\eta_\eps)|^2 + \frac1{\eps^2} V \eta_\eps^2 v^2 &=&  \int v^2 \eta_\eps  \Big(- \Delta \eta_\eps + \frac1{\eps^2} V \eta_\eps \Big) + \eta_\eps^2 |\nabla v|^2 \nonumber \\
	 &=& \int v^2 \eta_\eps  \Big(- \Delta \eta_\eps + \frac1{\eps^2} V \eta_\eps +  \frac1{\eps^2} \eta_\eps^3 \Big) -  \frac1{\eps^2} \eta_\eps^4 v^2 + \eta_\eps^2 |\nabla v|^2 \nonumber \\
	 &=& \frac {\lambda_\eps}{\eps^2}   \int v^2 \eta_\eps^2 +  \int  \eta_\eps^2 |\nabla v|^2 - \frac1{\eps^2} \eta_\eps^4 v^2  \label{eqenergytemp} \\
	 &=& \frac {\lambda_\eps}{\eps^2}  +  \int \eta_\eps^2 |\nabla v|^2  - \frac1{\eps^2} \eta_\eps^4 v^2 \nonumber \,.
 \ea

Using again (\ref{eqetaeps}), together with the mass constraint for $\eta_\eps$, we have that

\be \label{eqmueps}
	 \frac {\lambda_\eps}{\eps^2}  = 2 \Big( E_\eps(\eta_\eps) + \frac 1{4\eps^2} \int  \eta_\eps^4  \Big) \,.
\ee

Replacing (\ref{eqmueps}) in (\ref{eqenergytemp}), and then (\ref{eqenergytemp}) in (\ref{Eetaepsv}) we get

\ban
	 \E_{\eps}(u_1, u_2)  &=&  E_\eps(\eta_\eps) + \frac12 \int \eta_\eps^2 |\nabla v|^2 +   \frac 1{2\eps^2}  \eta_\eps^4 \{1 -2v^2 \} \\
	 &&   + \frac 12  \int  \frac 1{4}  v^2\eta_\eps^2  |\nabla  \varphi |^2 + \frac1{4\eps^2}  \eta_\eps^4 v^4 \{ 1 + \cos^2 \varphi \} + \frac 1{4}   g_\eps  \, \eta_\eps^4 v^4 \{ 1 - \cos^2 \varphi \} \,.
 \ean

Completing the square for $\{ 1 - v^2\}$ we get

\ban
	 \E_{\eps}(u_1, u_2)  &=&  E_\eps(\eta) + \frac12 \int \eta_\eps^2 |\nabla v|^2 + \frac1{2\eps^2}  \eta_\eps^4 \{1-v^2 \}^2 -  \frac1{2\eps^2}  \eta_\eps^4v^4  \\
	 &&   + \frac 12  \int   \frac 1{4}  v^2\eta_\eps^2  |\nabla  \varphi |^2 +  \frac1{4\eps^2}  \eta_\eps^4 v^4 \{ 1 + \cos^2 \varphi \} + \frac 1{4}   g_\eps  \, \eta_\eps^4 v^4 \{ 1 - \cos^2 \varphi \} \\
	       			              &=& E_\eps(\eta) + \frac12 \int \eta_\eps^2 |\nabla v|^2 + \frac1{2\eps^2}  \eta_\eps^4 \{1-v^2 \}^2  \\
	 &&   + \frac 12  \int   \frac 1{4}  v^2\eta_\eps^2  |\nabla  \varphi |^2 +  \frac 1{4}  \eta_\eps^4 v^4  g_\eps \Big( 1 - \frac1{g_\eps\eps^2} \Big)  \, \{ 1 - \cos^2 \varphi \} \,,
 \ean

which finishes the proof.\\  \qed

%%%%%%%%%%%%%%%%
%%%%%%%%%%%%%%%%
\section{Upper bound inequality } \label{upperbound}

In this section, we consider the formulation of the problem in $(v,\vp)$ and call
$$\F_\eps (v,\vp)=F_\eps(v)+G_\eps (v,\vp).$$
We prove here the upper bound inequality for $\eps \F_\eps$:

\begin{prop} \label{propUB}

\textbf{(Upper bound inequality for $\eps \F_\eps$)} Let $\vp=\pi \mathbf{1}_A \in X$. There is a sequence of pairs $(v_\eps, \vpeps) \in Lip(\R^2;(0,1] \times [0,\pi])$, converging as $\eps \to 0$ to $(1,\vp) $ in $L^1_{\text{loc}}(\D) \times L^1_{\text{loc}}(\D)$, such that
\ben
	\limsup_{\eps \to 0} \eps \F_\eps (v_\eps, \vpeps) \leq \F(\vp)  \,.
\een
\end{prop}

The proof is based on Bouchité's paper \cite{BOU}, where he proves the $\G$-convergence of an anisotropic phase transition Cahn-Hilliard energy.  We point out that our weight $\eta_\eps$ depends on $\eps$ and vanishes asymptotically on the boundary of $\cal D$.

In a first step, we assume that  $\vp = \pi \mathbf{1}_A$, where $A$ is an open bounded subset of $\R^2$ with smooth boundary such that $\H^1(\pA \cap \pD)=0$. Then, for any $\vp\in X$ we approximate $\{\vp=\pi \}$ by this kind of sets. We conclude then thanks to a density argument. We remark that we do not consider here the mass constraints in (\ref{mass2}).\\

Before proving the upper bound, we recall some results about sets with smooth boundary, that can be found in Lemmas 3 and 4 of \cite{Mod}. For an open set $A \subset \R^2$ with smooth, non empty compact boundary, let $d$ be the signed distance to $\pA$, defined by

\ben
	d(x)=  \left\{ \begin{array}{rll}
		 \text{dist}(x,\pA)& \text{ if }& x \in A \\
		 -\text{dist}(x,\pA) & \text{ if }&x \in \R^2\backslash A \,.
\end{array} \right.
\een

For small $t>0$, consider the neighborhood of $\pA$ given by

\ben
	N_t = \{ x \in \R^2 \,;\, |d(x)|<t \} \,,
\een

with boundary

\ben
	S_t = \{ x \in \R^2 \,;\, |d(x)|=t \} \,.
\een

For $t>0$ small enough, there is a diffeomorphism $\Phi$ between $N_t$ and $\pA \times ]0,t[$ such that

\be \label{bDPhi}
	\exists b >0 \,, \quad \det{|D\Phi|} \geq b \,.
\ee

We denote by $\hat \Phi$ the component of $\Phi$ in $\pA$. Moreover, $d$ is a Lipschitz continuous function in $N_t$ and we have that

\be \label{distderiv}
	|Dd| = 1 \quad \text{ a.e.\ in }  N_t \,.
\ee

For small $t>0$, define the measure

\ben
	\mu_t = \H^1 \mesrest (\D \cap S_t) \,.
\een

Notice that $\mu_0 =  \H^1 \mesrest (\D \cap \pA)$. As in Lemma 4 in \cite{Mod}, $\H^1(\partial A\cap \pD)=0$ yields

\ben
	\liminf_{t\to0} \mu_t(\Omega) \geq \mu_0(\Omega) \,,
\een

for every open $\Omega \subset \R^2$, and

\be \label{limmuD}
	\lim_{t\to0} \mu_t(\D) = \mu_0(\D) \,.
\ee

Hence, as $t \to 0$, $\mu_t$ converges weakly$^*$ to $\mu_0$, which implies

\be \label{limsupmeasure}
	\limsup_{t \to 0} \int_\D u \,  d\mu_t \leq  \int_\D u \,  d\mu_0  \,
\ee

for every upper semicontinuous function $u : \D \to \R$ with compact support (see Propositions 1.62 and 1.80 in \cite{AFP}). \\

Denote $\eta_0=\sqrt \rho$ and for $\eps\geq0$ define $f_\eps:\R^2\times\R\times\R^2\to\R_+$ by

\ben
	f_\eps(x,t,p)=  \frac12\eta_\eps^2(x) |p|^2 +  \frac1{4} \eta_\eps^4(x) \{1-t^2\}^2 \,.
\een

For $|p|=1$ and $s \in \R$ we also write $f_\eps(x,t,s)=f_\eps(x,t,sp)$.\\

The last step in the proof of  Proposition \ref{propUB} uses the following Lemma, which proof is given in the appendix.

\begin{lemm} \label{approxA}
	Let $A$ be a subset of $\D$ with $\mathbf{1}_A \in BV_{\text{loc}}(\D)$. There exists a sequence $\{A_k\}_{k \in \N}$ of open bounded subsets of $\, \R^2$ with smooth boundaries such that:
	\begin{description}
	\item[(i)] $\lim_{k \to \infty} \L^2((A_k \cap\D)\Delta A)=0$,
	\item[(ii)]  $\limsup_{k \to \infty} \int_{\D} \rho^{\nicefrac32} \, |D\mathbf{1}_{A_k}| \leq \int_{\D} \rho^{\nicefrac32} \, |D\mathbf{1}_A|$,
	\item[(ii)]   $\int_{A_k \cap \D} \rho  = \int_{A} \rho$ and $\H^1(\pD \cap \partial A_k)=0$ for $k$ large enough.
\end{description}
\end{lemm}

\textbf{Proof of Proposition \ref{propUB}:} we first assume that $A$ is an open subset of $\R^2$ with smooth, non empty compact boundary such that

\be \label{papD0}
	\H^1(\partial A\cap \pD)=0 \,.
\ee

 \textit{(Step 1: construction of the pairs of test functions.)}  For $T>1$, consider the approximation of the optimal profile

\ben
	w_T  = \left\{ \begin{array}{cll}
		\tanh &   \text{ in } &  (0,  T) \\
		h &   \text{ in } &   (T,T+\nicefrac1T)  \qquad \,,  \\ 		
		1 &   \text{ in }  & (T+\nicefrac1T,+\infty )
\end{array} \right.
\een

where $h$ is the unique cubic polynomial such that $h(T)=\tanh(T)$,  $h'(T)=\tanh'(T)$, $h(T+\nicefrac1T)=1$ and  $h'(T+\nicefrac1T)=0$. Computing explicitly the coefficients of $h$, we find that $w_T$ is a nondecreasing function in $\R_+$, with uniform $C^1$-bounds with respect to $T \in (1,\infty)$. We extend $w_T$ to the whole real line by setting $w_T(t) = w_T(-t)$ in $\R_-$. For $\eps \geq 0$, consider $0<m_\eps \ll \eps$ (to be chosen later), $t_\eps = \tanh^{-1}(m_\eps)$ and define a modification of $w_T$ near zero by

\ben
	w_{\eps,T}  = \left\{ \begin{array}{cll}
			m_\eps  &   \text{ in } &  (0, t_\eps)  \\
			w_T  &   \text{ in } &  (t_\eps,\infty)\, .
\end{array} \right.
\een

Notice that $w_{\eps,T}$ has uniform Lipchitz bounds with respect to $T \in (1,\infty)$ and $\eps \in [0,1)$. We recall that $\D=B(0,\lambda)$ and we denote $ \D_\delta = \overline{B(0,\lambda-\delta)}$. For $y$ in $\pA \backslash \D_\delta$, we define

\ben
	w_{\eps,T}^y(t) = w_{\eps,T}\left( \sqrt{\frac {\rho(y)}2} \, t\right) \,,
\een

and we write $w_T^y=w_{0,T}^y$. For small $\delta>0$, we define $R = R_\delta$ by

\be \label{defRdelta}
	R= \left( T +\nicefrac1T\right){\sqrt{\frac{2}{\delta \lambda}}} \,.
\ee

Since $w_T$ has uniform $C^1$-bounds with respect to $T \in (1,\infty)$, while $\rho$ is a smooth function in $\D_\delta$, for every  $y \in \pA \cap \D_\delta$, there is an open neighborhood $\S$ of $y$ in $\pA  \cap  \D_\delta$  such that %\footnote{We can take for example $V =B(y,\delta^2) \cap \D_\delta$.}

\ben
	\int_0^{R} f_0(x,w^y_T(t),(w^y_T)'(t)) \, dt \leq \int_0^{R}  f_0(x,w^x_T(t),(w^x_T)'(t)) \, dt \, +\,  \delta  \qquad \forall x \in \S \,, \quad  \forall T \geq 1 \,.
\een

Hence, thanks to the compactness of $\pA \cap \D_\delta$, there is a finite family $\{ \S_i\}_{i=1}^N$ of open disjoint subsets of $\pA \cap \D_\delta$, and a corresponding family of points $y_i \in \S_i$, such that

\be \label{H1meausreSi}
	\H^1\left( \pA \cap \D_\delta \setminus \bigcup_{i=1}^N \S_i \right) = 0 \,
\ee

and

\be \label{estf0Si}
	\int_0^{R} f_0(x,w^{y_i}_T(t),(w^{y_i}_T)'(t)) \, dt \leq \int_0^{R} f_0(x,w^x_T(t),(w^x_T)'(t)) \, dt \, +\,  \delta   \,,
\ee

for every $x \in \S_i$,  $T \geq 1$ and $1\leq i \leq N$. We will use the functions $w^{y_i}_{\eps,T}$ to define the first test function, so we have to interpolate between the different $\S_i$'s. Define first  $\S_0 = \pA \backslash \D_\delta$ and $y_0 = (\lambda-\delta,0) \in \partial \D_\delta$. For small $\l>0$ define $\S_i^\l = \{ x \in \S_i \,;\, \text{dist}(x,\partial\S_i)  \geq \l\}$. Clearly,

\ben
	\H^1(\S_i \backslash \S_i^\l) \to 0 \quad \text{ as} \quad \l \to 0.
\een

In particular, we can take $\l = \l_\delta$ such that

\be \label{ldelta}
	R \, \H^1(\S_i \backslash \S_i^\l) = o_{\delta \to 0}(1) \,
\ee

for every $0\leq i \leq N$. Consider then $\{ \hat \theta_i \}_{i=0}^N$ such that $\hat \theta_i\in C^\infty(\pA,[0,1])$,

\be \label{interpolation}
	\sum_{i =0}^N\hat \theta_i =1 \quad \text{on } \quad \pA  \qquad \text{ and } \qquad \hat \theta_i=1 \quad \text{in } \quad \S_i^\l \,.
\ee

We deduce a smooth partition of the unity on $N_t$ by setting $ \theta_i =\hat \theta_i \circ \hat \Phi$ and we define

\ben
	v_\eps =  \left\{ \begin{array}{ccl}
			1 & \text{ in } &  \R^2 \backslash N_{\eps R}  \\ \\
			\sum_{i =0}^N \theta_i (x) \, w^{y_i}_{\eps,T} \left( \frac{|d(x)|}\eps \right) & \text{ in } &  N_{\eps R}
	\end{array}\right. \,.
\een

Since $w^{y_i}_{\eps,T}$ is a nondecreasing function, while $\rho$ is a radial decreasing function,  (\ref{estrhopD}) and the fact that $\text{dist}(y_i,\D) \geq \delta$ yield

\be \label{contveps}
	w^{y_i}_{\eps,T} \big|_{\partial N_{\eps R}}  = w^{y_i}_{\eps,T} (R) = w_{\eps,T} \left(    \left( T +\nicefrac1T\right){\sqrt{\frac{\rho(y_i)}{\delta \lambda}}}  \right) \geq  w_{\eps,T} (  T +\nicefrac1T )= 1 \,,
\ee

so $v_\eps$ is a continuous function. Moreover, since  $w_{\eps,T}$ has uniform Lipschitz bounds with respect to $T \in (1,\infty)$ and $\eps \in (0,1)$, there is $C>0$ such that for $\eps $ small enough,

\be \label{estdervy}
	\| v_\eps \|_{C^{0,1}(\R^2)} \leq   \frac C\eps  \,.
\ee

We also define

\be \label{defphieps}
	\vp_\eps(x) = \left\{ \begin{array}{ccl}
	\pi & \text{ if } &  x \in  A \backslash N_{\eps \tilde t_\eps} \\
	  \xi(\nicefrac{d(x)}{\eps t_\eps}) & \text{ if } & x \in   N_{\eps \tilde t_\eps} \\
	0 & \text{ if } & x  \in  \R^2  \backslash ( A \cup N_{\eps \tilde t_\eps})
	\end{array}\right. \,,
\ee

where $\xi(t)=\nicefrac\pi2(1+t)$ and $\tilde t_\eps =  (\nicefrac 2 \lambda)^{\nicefrac12} \, t_\eps$. We clearly have that $(v_\eps, \vpeps) \in Lip(\R^2;(0,1] \times [0,\pi])$, and that $(v_\eps, \vpeps)$ converges as $\eps \to 0$ to $(1,\vp) $ in $L^1_{\text{loc}}(\D) \times L^1_{\text{loc}}(\D)$.\\

%%%%%%
\textit{(Step 2: estimating the energy $\eps G_\eps$.)} The function $\vp_\eps$ is constant in $\R^2\backslash N_{\eps \tilde t_\eps}$, so $G_\eps(v_\eps,\vp_\eps)  = G_\eps(v_\eps,\vp_\eps \, ;\, N_{\eps \tilde t_\eps} ) $. Since $w_{\eps,T}$ is a nondecreasing function while $\rho$ has a global maximum at zero, for every  $x \in  N_{\eps \tilde t_\eps}$ we have

\ban
	w^{y_i}_{\eps,T} \left( \frac{|d(x)|}\eps \right) \leq w^0_{\eps,T} \left(\frac{|d(x)|}\eps\right) = w_{\eps,T} \left( \frac \lambda  {\sqrt{2}} \frac{|d(x)|}\eps \right) \leq  w_{\eps,T} ( t_\eps) = m_\eps \,,
\ean

so $v_\eps\leq m_\eps $ in $N_{\eps \tilde t_\eps}$. Hence,

\ben
	 G_\eps(v_\eps,\vp_\eps) \leq  \frac18 \int_{N_{\eps \tilde t_\eps}}  \eta_\eps^2 m_\eps^2  \, |\nabla \varphi_\eps|^2 + \eta_\eps^4 m_\eps^4 \,  \tilde g_\eps  \{ 1-\cos^2(\varphi_\eps)\} \,.
\een

Then, the definitions of $\vp_\eps$ and $ \tilde g_\eps$, together with  the fact that $\eta_\eps$ is uniformly bounded, yield
\ben
	  G_\eps(v_\eps,\vp_\eps) \leq
	C   \left( m_\eps^2 \, (\eps \tilde t_\eps)^{-2} +  g_\eps \, m_\eps^4  \right) \, \L^2(N_{\eps \tilde t_\eps}) \,.
\een

Using (\ref{bDPhi}), we have

\ban
	(\eps \tilde t_\eps)^{-1} \L^2(N_{\eps \tilde t_\eps})  &\leq & C \,(\eps \tilde t_\eps)^{-1}  \,  \int_{-\eps \tilde t_\eps}^{ \eps \tilde t_\eps}  \, dt \int_{\pA }  \left| \det{(D\Phi)^{-1}}  \right| \, d\H^1 \\
	 &\leq&     b^{-1} \, \H^1 (\pA)  \\
	 &=& \O_{\eps \to 0}(1) \,.
\ean

For $\eps$ small, we have $ \tilde t_\eps = ( \nicefrac   2 \lambda)^{\nicefrac12} \tanh^{-1}(m_\eps) =  \O_{\eps \to 0}(m_\eps)$. Hence,

\ben
	  G_\eps(v_\eps,\vp_\eps) \leq	C  \left( m_\eps \eps^{-1}+ m_\eps^5 \,  g_\eps \, \eps \right)  \,.
\een

Taking $m_\eps = (g_\eps \eps^2)^{-\nicefrac14}  $, $G_\eps(v_\eps,\vp_\eps)  \leq	C  g_\eps^{-\nicefrac14}  \eps^{-\nicefrac32}$, and after (\ref{regime}) we obtain

\be \label{limepsGeps}
	\eps G_\eps(v_\eps,\vp_\eps) \leq C  (g_\eps \eps^2)^{-\nicefrac14}    = o_{\eps\to0}(1) \,.
\ee

%%%%%
\textit{(Step 3: computing the energy $\eps F_\eps$.)} Since $v_\eps$  is constant out of $N_{\eps R} $, we have

\be \label{energyVeps}
	\eps F_\eps(v_\eps) = \int_{N_{\eps R} \cap \D_\delta }   \phi_\eps  + \int_{N_{\eps R} \cap ( \D \backslash \D_\delta) }  \phi_\eps  +\int_{N_{\eps R} \backslash \D }  \phi_\eps \,.
\ee

where

 \ben
	\phi_\eps(x) = \frac1\eps f_\eps(x,v_\eps,\eps Dv_\eps) \,.
\een

Considering (\ref{estdervy}), for $\eps$ small enough there is $C>0$ such that $| \nabla v_\eps | \leq  \nicefrac C\eps $. Estimates (\ref{estetaK}), (\ref{estetarad}) and (\ref{estrhopD}) thus yield

\ban
	\int_{N_{\eps R} \cap ( \D \backslash \D_\delta)} \, \phi_\eps  &\leq& C \, \Big( \max_{\partial \D_\delta} \{\rho+ \rho^2\} + C_\delta \, \eps^2 |\ln \eps|\Big) \, \eps^{-1} \, \L^2(N_{\eps  R } \cap ( \D \backslash \D_\delta) ) \\
	&\leq& C  \, \big( \delta + \eps^2 |\ln \eps| \big) \, \eps^{-1}  \, \L^2(N_{\eps R} ) \,.
\ean

Using (\ref{bDPhi}) we have

\ban
	\eps^{-1} \, \L^2(N_{\eps R} )  &\leq & C \, \eps^{-1} \int_{-\eps R}^{ \eps R}  \, dt \int_{\pA }  \left| \det{(D\Phi)^{-1}}  \right| \, d\H^1 \\
	 &\leq&   C \, R \, b^{-1} \, \H^1 (\pA)   \,,
\ean

so (\ref{defRdelta}) yields

\be  \label{energyVeps2}
	\lim_{\eps \to 0} \int_{N_{\eps R} \cap ( \D \backslash \D_\delta)} \, \phi_\eps \leq  C  \,  \delta  \,  R \, b^{-1} \, \H^1 (\pA) =o_{\delta \to 0}(1)  \,.
\ee

Similarly, using (\ref{estetaoutbulk}) we have

\ba
	\int_{N_{\eps R}  \backslash \D} \, \phi_\eps  &\leq&  C  \, \sup_{\R^2 \backslash \D}  \big( \eta_\eps^2 + \eta_\eps^4  \big)\,  R \, b^{-1} \, \H^1 (\pA)  \nonumber \\
	 									&\leq& C \,  \eps^{\nicefrac13}  \,  R \, b^{-1} \, \H^1 (\pA)  \label{energyVeps12} \\
									     &=& o_{\eps \to 0}(1)\,. \nonumber
\ea

Now, remember the interpolation from (\ref{interpolation}). We have

\ben
	\int_{N_{\eps R} \cap \D_\delta } \, \phi_\eps  =  \sum_{i =1}^N \left(  \int_{N_{\eps R} \cap B_i }  \phi_\eps +  \ \int_{N_{\eps R} \cap C_i }  \phi_\eps\right)  \,,
\een

where

\ben
	B_i=\{ x \in \D_\delta \,;\, \hat \Phi (x) \in \S_i^\l \} \qquad \text{ and } \qquad C_i= \{x \in \D_\delta \,;\, \hat \Phi  (x)\in \Sigma_i \backslash \S_i^\l \} \,.
\een

As before, since $\eta_\eps$ is uniformly bounded in $\R^2$ with respect to $\eps \in(0,1)$,  (\ref{bDPhi}) yields

\ban
	 \int_{N_{\eps R} \cap C_i }\phi_\eps &\leq & C \, \frac1\eps \int_{-\eps R}^{ \eps R}  \, dt \int_{ \Sigma_i \backslash \S_i^\l}  \left| \det{(D\Phi)^{-1}}  \right| \, d\H^1 \\
	 &\leq&   C \, R \, b^{-1} \, \H^1 (\Sigma_i \backslash \S_i^\l)   \,.
\ean

Hence, after (\ref{defRdelta}) and  (\ref{ldelta}), we obtain

\be  \label{energyCi}
	\int_{N_{\eps R} \cap C_i }\phi_\eps = o_{\delta \to 0}(1)  \,
\ee

for every $\eps \in (0,1)$. In $N_{\eps R} \cap B_i$ we have $ v_\eps(x)=w^{y_i}_{\eps,T} (\nicefrac{|d(x)|}\eps)$.  Using (\ref{distderiv}) we write

\ben
	\int_{N_{\eps R} \cap B_i} \phi_\eps  = \int_{N_{\eps R}  \cap  B_i} |D(\nicefrac{|d|}\eps)| (x)  f_\eps(x,w^{y_i}_{\eps,T} (\nicefrac{|d(x)|}\eps),(w^{y_i}_{\eps,T})'  (\nicefrac{|d(x)|}\eps)) \,.
\een

The coarea formula from Proposition \ref{coareaprop} yields

\ben
	\int_{N_{\eps R} \cap B_i} \phi_\eps   =   \int_{-R}^{R}  dt \int_{\D_\delta \cap B_i} f_\eps(x,w^{y_i}_{\eps,T} (t),(w^{y_i}_{\eps,T})' (t))\,  d \mu_{\eps t}(x) \,.
\een

We thus have,

\be \label{errorRi}
	\int_{N_{\eps R} \cap B_i} \phi_\eps \leq \int_{-R}^{R} dt \int_{\D_\delta \cap B_i} f_0(x,w^{y_i}_{T} (t),(w^{y_i}_{T})' (t))\,  d \mu_{\eps t}(x)  + R^i_{\eps} + \tilde R^i_{\eps}   \,.
\ee

The first error here before comes from the modification of $w_T$ near $0$. Using (\ref{limmuD}) and the definition of  $t_\eps$ we compute

\ban
	 R^i_{\eps} &\leq& C  \, \sqrt{\frac2{\rho(y_i)}} \, t_\eps    \, \sup_{t \in(0,\eps R) } \| \mu_{ t } \|   \\
	 		    &\leq& C  \, \sqrt{\frac2{\delta \lambda}} \, t_\eps    \, \sup_{t \in(0,\eps R) } \| \mu_{ t } \|   \\
			     &=& o_{\eps \to 0}(1) \,.
\ean

The second error appears when replacing $f_\eps$ by $f_0$, so using estimates (\ref{estetaK}) and  (\ref{estrhopD}), together with $y_i \in \D_\delta$, there is $C_\delta >0 $ such that

\ban
	 R^2_{\eps} &\leq C_\delta \,  \eps^2 \, |\ln \eps| \, R   \,   \, \sup_{t \in(0,\eps R ) } \| \mu_{ t } \|  =  o_{\eps \to 0}(1) \,.
\ean

Using Fubini's formula, we rewrite (\ref{errorRi}) as

\ban
	\int_{N_{\eps R} \cap B_i} \phi_\eps  &\leq&  \int_{\D}  \left( \mathbf{1}_{\D_\delta \cap B_i}(x) \int_{-R}^{R}  f_0(x,w^{y_i}_T (t),(w^{y_i}_T)' (t) )\,  dt \right) d \mu_{\eps t}(x)  + o_{\eps \to 0}(1) \,.
\ean

The set $ \D_\delta \cap B_i $ is close and the inner integral is a continuous function of $x$. Hence, the function inside the outer integral is upper semicontinuous function of $x$. Inequality (\ref{limsupmeasure}) thus yields

\ban
	\limsup_{\eps \to 0} \int_{N_{\eps R} \cap B_i} \phi_\eps &\leq&   \int_{\D}   \left( \mathbf{1}_{\D_\delta \cap B_i}(x) \int_{-R}^{R}  f_0(x,w^{y_i}_T (t),(w^{y_i}_T)' (t) )\,  dt \right) d \mu_{0}(x)   \\
												&=&  \int_{\D_\delta \cap \S_i^\l}    \left(  \int_{-R}^{R}  f_0(x,w^{y_i}_T (t),(w^{y_i}_T)' (t) )\,  dt \right) d \mu_{0 }(x) \,.
\ean

Notice that since $\mu_0$ is supported in $\pA$, we replaced $B_i$ by $\D_\delta \cap \S_i^\l$. From (\ref{estf0Si}) and since $w_T$ is an even function, we have

\ba
	\limsup_{\eps \to 0} \int_{N_{\eps R} \cap B_i} \phi_\eps &\leq&   2   \int_{\D_\delta \cap \S_i^\l}    \left(  \int_0^{R}  f_0(x,w^{x}_T (t),(w^{x}_T)' (t) )\,  dt \,+\, \delta \right) d \mu_{0 }(x) \nonumber \\								 &\leq&  2   \int_{\D_\delta \cap \S_i^\l}    \left(  \int_0^{\infty}  f_0(x,w^{x}_T (t),(w^{x}_T)' (t) )\,  dt \,+\, \delta \right) d \mu_{0 }(x)  \,. \label{energyVeps0}
\ea

\textit{(Step 4: upper bound)} Putting together (\ref{limepsGeps})-(\ref{energyCi}) and (\ref{energyVeps0}) we get

\ban
	 \limsup_{\eps \to 0}\eps \F_\eps(v_\eps,\vp_\eps)  &\leq&   2 \sum_{i =1}^N   \int_{\D_\delta \cap \S_i^\l}    \left(  \int_0^{\infty}  f_0(x,w^{x}_T (t),(w^{x}_T)' (t) )\,  dt \,+\, \delta \right) d \mu_{0 }(x)    + o_{\delta \to 0}(1) 	  \,.
\ean

Now, we take a sequence $T= T_\delta$ such that $T_\delta \to \infty$ as $\delta \to 0$ (notice that (\ref{contveps}) still holds). Then, Fubini's formula and dominated convergence theorem, together with (\ref{H1meausreSi}) and (\ref{ldelta}), yield

\ban
	\lim_{\delta \to 0} \left( \limsup_{\eps \to 0}\eps \F_\eps(v_\eps,\vp_\eps) \right) &\leq&     \int_{\D}  \left( \int_0^\infty  f_0(x,w^{x} (t),(w^{x})' (t) )\,  dt  \right)  d \mu_{0}(x)\,.
\ean

Remembering the definitions of $f_0$, $\mu_0$ and $w^{x}$, (\ref{infI}) yields

\ban
	\lim_{\delta \to 0} \left( \limsup_{\eps \to 0}\eps \F_\eps(v_\eps,\vp_\eps)  \right) &\leq&    2 \sigma   \int_{\D \cap \pA}    \rho^{\nicefrac32}(x)   d \H^1(x)	  \,.
\ean

We conclude thanks to a diagonal argument (see Corollary 1.16 in \cite{attouch}): there exists a sequence $\delta_\eps \to_{\eps \to 0} 0$, such that as $\eps \to 0$, $(v_{\eps,\delta_\eps},\vp_{\eps,\delta_\eps})$ converges in $L^1_{\text{loc}}(\D) \times L^1_{\text{loc}}(\D)$ to $(1,\vp)$, and

\ban
	\limsup_{\eps \to 0}\eps \F_{\eps,\delta_\eps}(v_{\eps,\delta_\eps},\vp_{\eps,\delta_\eps})   &\leq&    2 \sigma   \int_{\D \cap \pA}    \rho^{\nicefrac32}(x)   d \H^1(x)	  \,.
\ean

%%%%%%%%
\textit{(Step 5: approximation of $A$ by Cacciopoli sets)} We end the proof using Lemma \ref{approxA}, the proof of which is given in the appendix. We remove the condition (\ref{papD0}) and we only assume that $A$ is a set with locally finite perimeter in $\D$. Consider $\vp^k=\pi \mathbf{1}_{A_k}$ with $\{A_k\}_{k \in \N}$ the sequence from Lemma \ref{approxA}. From \textbf{(i)}, $(1,\vp^k)$ converges to $(1,\vp)$ in $L^1(\D)$ and we have $\int_{A_k} \rho= \alpha_2$ for $k$ large enough. Hence, from steps (1)-(4), there is a sequence $(v_\eps^k,\vp_\eps^k) \to (1,\vp^k)$ in $L_{\text{loc}}^1(\D) \times L_{\text{loc}}^1(\D)$ such that

 \ben
	\lim_{\eps \to 0} \eps \F_\eps (v_\eps^k,\vp_\eps^k) = \F(\vp^k) \,.
\een

Using \textbf{(ii)} from Lemma \ref{approxA} we obtain

 \ben
	\limsup_{k \to 0} \left( \lim_{\eps \to 0} \eps \F_\eps (v_\eps^k,\vp_\eps^k) \right) \leq \F(\vp) \,.
\een

As in step (4), we conclude thanks to a diagonal argument. \\ \qed \\

%%%%%%%%%%%%%%%%
%%%%%%%%%%%%%%%%
\section{Lower bound inequality and compactness} \label{lowerbound}

The proofs in this section are based on geometric measure theory techniques. We make the lower bound on lines and then use the slicing method, which can be found in \cite{AT} or \cite{braides}. The last part of the proof of the lower bound is inspired by the ideas in \cite{ABS}.

%%%%%%%
\subsection{Lower bound on lines}

Consider an open set $A \subset \subset \D$ and let $\nu \in S^1$ be a fixed direction. We call $\pi_\nu$ the hyperplane orthogonal to $\nu$, and $A_\nu$ the projection of $A$ on $\pi_\nu$. We define the one dimensional slices of $A$, indexed by $x \in A_\nu$, as

\ben
	 A_x = \{ t \in \R \,;\, x+t\nu \in A\} \,.
\een

 For every function $f$ in $\D$, we define $f_x$ as the restriction of $f$ to the slice $A_x$, defined by $f_x(t) = f(x+t\nu)$. For $(v,\vp) : A_x \to (0,1] \times (0,\pi)$, we define the energies

\ban
	F_{\eps}(v \,;A_x) &=&   \frac12  \int_{A_x} \eta_{\eps,x}^2 |\nabla v|^2 +  \frac1{2\eps^2} \eta_{\eps,x}^4 \{1-v^2\}^2 \,, \\
	G_\eps(v,\varphi \,;A_x) &=&  \frac18   \int_{A_x}    \eta_{\eps,x}^2 v^2  \, |\nabla \varphi|^2 + \eta_{\eps,x}^4 v^4 \, \tilde g_\eps  \{ 1-\cos^2(\varphi)\} \quad \text{ and } \\
	\F_\eps(v,\varphi\,;A_x) &=& F_\eps(v;A_x) + G_\eps(v,\varphi;A_x) \,.
\ean

Similarly, for $\vp  \in BV(A_x\,;\{0,\pi\})  $ we define

\ben
	\F(\vp \, ; A_x) = \frac{2\sigma}{\pi} \, \int_{A_x} \rho_x^{\nicefrac32} d |D\vp|   \,.
\een

With the previous notations, we have the following result:
%%%%%%%%%
\begin{prop}\label{lowerboundonlines} Let $(v_\eps,\vp_\eps) \in Lip(A_x \,;\, (0,+\infty) \times [0,\pi])$ such that

\be \label{boundlinftv}
	\sup_{\eps>0} \|v_\eps\|_{L^{\infty}(A_x)} < \infty
\ee

and

\be \label{unifbound1d}
	\eps \F_\eps(v_\eps,\vp_\eps \, ;A_x)< \infty \,.
\ee

Then, there is $\vp \in SBV(A_x\,;\{0,\pi\})$ such that

\ba
	(v_\eps,\vp_\eps) \to (1,\vp) \quad &\text{ in }& \quad L^1(A_x) \times L^1(A_x) \label{lw1d1}
\ea

and

\be
	\liminf_{\eps\to0} \eps \F_\eps(v_\eps,\vp_\eps\, ;A_x) \geq \F(\varphi \, ;A_x) \,. \label{lw1d3}
\ee

\end{prop}

\textbf{Proof:} \textit{(Step 1)} Using that $A \subset\subset\D $ and estimate (\ref{estetaK}), there are $c_1,c_2>0$ such that

\be \label{lbetaepsrho1d}
	 \eta_{\eps,x}^2 > \rho_x - c_1\eps^2 |\ln \eps|  > c_2 \quad \text{ in } \quad A_x \,.
\ee

Hence, the definition of $\F_\eps( \cdot\,; A_x)$  and (\ref{unifbound1d}) give

\ben
	 \int_{A_x} |1-v_\eps|  < \frac{4 |A_x|}{c_2^2} \, \eps^2 \F_\eps(v_\eps,\vp_\eps \, ;A_x) = o_{\eps\to0}(1) \,,
\een

so $v_\eps \to 1$ in $L^1(A_x)$. Similarly, after (\ref{boundlinftv}) $v_\eps <C$ in $A_x$, so (\ref{regime}) yields

\ben
	 \int_{A_x} v_\eps^4 |1-\cos^2 (\vp_\eps) |  < \frac{8C}{\tilde g_\eps \, c_2^2} \, \F_\eps(v_\eps,\vp_\eps \, ;A_x) = o_{\eps\to0}(1) \,.
\een

Hence, up to a (not relabeled) subsequence,  $\vp_\eps \to \vp$ a.e.\ in $A_x$, with $\vp : A_x \to \{0,\pi\}$. This, together with $A_x \subset\subset \D$, gives $\vp_\eps \to \vp$ in $L^1(A_x)$. We have proved (\ref{lw1d1}).\\

\textit{(Step 2)}  We now prove the lower bound for the energy. Let $t_0 \in S\vp$. For $\delta>0$ define

\ben
	J_\delta = A_x \cap (t_0-\delta,t_0+\delta) \,,
\een

and suppose that

\be \label{t0deltalb}
	\inf_{t \in J_\delta} \left\{ \inf_{\eps >0} v_\eps(t) \right\} > c_3 >0 \,.
\ee

Then, for every $\eps>0$ and every $t \in J_\delta$, $v_\eps(t)>c_3$. Hence, using (\ref{unifbound1d}), (\ref{lbetaepsrho1d}) and  the coarea formula (\ref{coarea}), there is $C>0$ such that

\be  \label{ineqW}
	 (c_2^{-1}c_3^{-2} + c_2^{-2} c_3^{-4}) \frac {8\,C}{\tilde g_\eps \eps^2}  \geq \int_{J_\delta}  |\nabla \varphi|^2 +   \{ 1-\cos^2(\varphi)\} \geq \int_0^\pi dt \, \{1-\cos^2(t)\}^{\nicefrac12} \int_{J_\delta}  |D\mathbf{1}_{W_{\eps,t}}| \,,
\ee

where $W_{\eps,t} = \{ t \in A_x \,;\, \vp_\eps(x) <t\}$. Since $\vp_\eps$ converges to  $\vp$ a.e.\ in $A_x$, we get

\ben
	\mathbf{1}_{W_{\eps,t}} \to \mathbf{1}_{ \{ \vp = 0 \}} \quad \text{ in } \quad L^1(A_x) \,,
\een

for a.e.\ $t \in (0,\pi)$. Hence, the lower semicontinuity of the BV norm with respect to the $L^1$-convergence, together with (\ref{regime}), (\ref{ineqW}) and Fatou's lemma, gives

\ban
	0 &\geq& \int_0^\pi dt \, \{1-\cos^2(t)\}^{\nicefrac12}  \liminf_{\eps \to 0} \int_{J_\delta}  |D\mathbf{1}_{W_{\eps,t}}|  \\
	&\geq& \frac\pi2 \int_{J_\delta}  |D\mathbf{1}_{\{ \vp = 0 \}}| \,.
\ean

Thus,

\ben
	0 = \H^0(J_\delta \cap S\vp) \geq \H^0(\{ t_0\}) =1 \,.
\een

This contradiction implies that  (\ref{t0deltalb}) can not be satisfied. We derive that for every $\delta >0$, we may extract a subsequence (not relabeled), such that exists $\{ t_\eps\}_{\eps >0} \subset J_\delta$ with

\be \label{ssteps}
	t_\eps \to \tilde t_0 \in \overline{J_\delta} \qquad \text{ and } \qquad v_\eps(t_\eps) \to 0  \qquad \text{ as } \eps \to 0 \,.
\ee

 \textit{(Step 3)} For $\eps>0$, define $I^{\pm}_\eps = \{t \in J_\delta \,;\, \pm(t_\eps-t) <0 \}$ and $v_\eps^{\pm}:J_\delta \to (0,1]$ by

\ben
	v_\eps^{\pm}(t) = \mathbf{1}_{I^{\pm}_\eps} v_\eps(t_\eps) + \mathbf{1}_{I^{\mp}_\eps} v_\eps(t) \,.
\een

The definition of $\F_\eps$, estimate (\ref{lbetaepsrho1d}) and the fact that $v_\eps^{\pm}$ is constant in $I^{\pm}_\eps$ while equal to $v_\eps$ in $I^{\mp}_\eps$ yield

\ban
	 \sqrt2 \, \eps \F_\eps(v_\eps,\vp_\eps\, ;J_\delta) \geq    \int_{J_\delta}  \rho_{x}^{\nicefrac32} \Big( |(v^+_\eps)'| \, |1-(v^+_\eps)^2|   +  |(v^-_\eps)'| \, |1-(v^-_\eps)^2| \Big) + o_{\eps \to 0}(1)  \,.
\ean

Using the coarea formula (\ref{coarea}) we obtain
\ben
	 \sqrt2 \, \eps \F_\eps(v_\eps,\vp_\eps\, ;J_\delta)  \geq   \int_0^1  dt \, (1-t^2) \, \int_{J_\delta}  \rho_{x}^{\nicefrac32} \Big( |D\mathbf{1}_{V^+_{\eps,t}}| +  |D\mathbf{1}_{V^-_{\eps,t}}| \Big) + o_{\eps \to 0}(1)   \,,
\een

where $V^\pm_{\eps,t} = \{ t \in J_\delta \,;\, v^\pm_\eps <t\}$. Since $t_\eps \to \tilde t_0$, $v_\eps(t_\eps) \to 0$ and $v_\eps(t) \to 1$ a.e.\ in $J_\delta$, $\mathbf{1}_{V^\pm_{\eps,t} } \to \mathbf{1}_{I^\mp}$ in $L^1(J_\delta)$, where  $I^{\pm} = \{t \in J_\delta \,;\, \pm(\tilde t_0-t)\leq0 \}$. Hence, the lower semicontinuity of the BV norm with respect to the $L^1$-convergence  and Fatou's lemma give

\be \label{lbtildet0}
	  \liminf_{\eps \to 0}\eps \F_\eps(v_\eps,\vp_\eps\, ;J_\delta)  \geq   \frac1{\sqrt2}  \int_0^1  dt \, (1-t^2) \, \int_{J}  \rho_{x}^{\nicefrac32} \Big( |D\mathbf{1}_{I^-}| +  |D\mathbf{1}_{I^+}| \Big) =   \sigma \, 2\rho_{x}^{\nicefrac32}(\tilde t_0) \,.
\ee

Moreover, since $\rho_{x}^{\nicefrac32}\geq c_4>0$ in $A_x$, we have

\be \label{lbtildet02}
	  \liminf_{\eps \to 0}\eps \F_\eps(v_\eps,\vp_\eps\, ;J_\delta)  \geq   2 \sigma c_4 >0 \,.
\ee

\textit{(Step 4)} Let $\Gamma = \{t_1, \dots, t_n\}$, $n \in \N$, be any finite subset of $S\vp$. For $i  \in \{0,1,\dots,n\}$ we define

\ben
	J_\delta^i = A_x \cap (t_i-\delta,t_i+\delta) \,.
\een

Consider  $\delta'>0$ such that $J_\delta^i \cap J_\delta^j=  \emptyset$ for $ i\neq j$ and let $ \delta \in(0,\delta')$. From (\ref{lbtildet02}), we have

\ben
	\frac  {c_4}{2 \sigma }  \liminf_{\eps \to 0}\eps \F_\eps(v_\eps,\vp_\eps\, ;A_x)   \geq  n \,.
\een

Therefore, using (\ref{unifbound1d}) we derive that $n$ is bounded, so $S\vp$ is a finite set and $\vp \in SBV(A_x)$.\\

\textit{(Step 5)} Finally, write $S\vp= \{t_1, \cdots, t_N\}$, $N \in \N$. Reasoning as before, for $\delta'$ small enough and $\delta \in(0,\delta')$,  (\ref{lbtildet0}) gives

\ben
	  \liminf_{\eps \to 0}\eps \F_\eps(v_\eps,\vp_\eps\, ;A_x)  \geq  2\sigma  \sum_{i=1}^N  \rho_{x}^{\nicefrac32}(\tilde t_i) \,, \qquad \text{ with } \qquad \tilde t_i \in \overline{J_\delta^i} \,.
\een

Since $\rho_{x}$ is a continuous function, taking the limit $\delta \to 0$ in the previous inequality we obtain

\ben
	  \liminf_{\eps \to 0}\eps \F_\eps(v_\eps,\vp_\eps\, ;A_x)  \geq  2\sigma  \sum_{i=1}^N  \rho_{x}^{\nicefrac32}( t_i) = 2\sigma \int_{A_x} \rho_{x}^{\nicefrac32} \, d|D\vp|  \,.
\een

We have proved (\ref{lw1d3}). \\ \qed \\

%%%%%%%%%%
\subsection{The slicing method}

Using the slicing method, we now prove the compactness and the lower bound inequality for $\eps \F_\eps$.

\begin{prop} \label{propLB} \textbf{(Lower bound inequality and compactness for $\eps \F_\eps$)} Let $(v_\eps,\vp_\eps) \in Lip_{loc}( \R^2 \,;\, (0,+\infty) \times [0,\pi])$ such that

\be \label{boundlinftv2}
	\sup_{\eps>0} \|v_\eps\|_{L^{\infty}(K)} < C_K \qquad \text{ for } K \subset \subset \D
\ee
and
\be \label{unifboundF}
	\sup_{\eps >0 }\eps \F_\eps(v_\eps,\vp_\eps)< \infty \,.
\ee

Then, there is $\vp \in X$ such that

\ba \label{cvvphi}
	(v_\eps,\vp_\eps) \to (1,\vp) \quad &\text{ in }& \quad L^1_{loc}(\D) \times L^1_{loc}(\D) \label{lw1d2}
\ea

and

\be \label{lowerboundF}
	\liminf_{\eps\to0} \eps \F_\eps(v_\eps,\vp_\eps\, ) \geq \F(\varphi ) \,.
\ee

\end{prop}

\textbf{Proof:} arguing as in the proof of (\ref{lw1d1}) of Proposition \ref{lowerboundonlines}, there exists $\vp : \D \to \{0,\pi\}$ such that (\ref{cvvphi}) is satisfied. Consider an open set $A \subset\subset\D$, and fix $\nu \in S^1$. For $x \in A_\nu$, we define $(v_{\eps,x},\vp_{\eps,x}) : A_x \to (0,1] \times (0,\pi)$ by

\ben
	(v_{\eps,x},\vp_{\eps,x})(t) = (v_{\eps},\vp_{\eps})(x + t \nu) \,.
\een

Since $A \subset\subset\D$ and since $v_\eps$, $\eta_\eps$ are non vanishing continuous functions, for fixed $\eps>0$ (\ref{unifboundF}) yields

\ban
	C \geq  \eps \F_\eps(v_\eps,\vp_\eps\, ;A) &\geq&    \frac\eps2 \inf_{A} \eta_\eps^2 \int_A |\nabla v_\eps|^2+  \frac\eps8 \inf_{A} v_\eps^2\eta_\eps^2 \int_A |\nabla \vp_\eps|^2 \\
								&\geq& c_{\eps,A} \left\{  \int_A |\nabla v_\eps|^2 + \int_A   |\nabla \vp_\eps|^2 \right\}\,,
\ean

so $v_\eps$ and $\vp_\eps$ belong to $W^{1,2}(A)$. Hence (see \cite{EG}, Section 4.9.2),

\ben
	v_{\eps,x}'(t) = 	D_\nu v_{\eps}(x+t \nu) \qquad \text{ and } \qquad \vp_{\eps,x}'(t) = 	D_\nu \vp_{\eps}(x+t \nu)
\een

for a.e.\ $t \in A_x$, for $\L^1$- a.e.\ $x \in A_\nu$. Using then $|\nabla v_\eps|^2\geq|D_\nu v_\eps|^2$, we get the slicing inequality

\be \label{slicingineq}
	 \eps\F_\eps (v_\eps,\vp_\eps\,;A) \geq \int_{A_\nu} \eps \F_\eps (v_{\eps,x},\vp_{\eps,x}\,;A_x) \, dx \,.
\ee

 From (\ref{unifboundF}), for $\L^1$-a.e.\ $x \in A_\nu$, $\eps \F_\eps (v_{\eps,x},\vp_{\eps,x}\,;A_x) $ is uniformly bounded with respect to $\eps$. Thus, after Proposition \ref{lowerboundonlines}, for $\L^1$-a.e.\ $x \in A_\nu$ there is $\vp_x \in BV(A_x\,;\{0,\pi\})$ such that

\ba
	(v_{\eps,x},\vp_{\eps,x}) \to (1,\vp_x) \quad &\text{ in }& \quad L^1(A_x) \times L^1(A_x) \label{lw1d1x}
\ea

 and

\be
	\liminf_{\eps>0} \eps \F_\eps(v_{\eps_x},\vp_{\eps,x}\, ;A_x) > \F(\varphi_x \, ;A_x) \,. \label{lw1d3x}
\ee

The function $\vp$ defined in (\ref{cvvphi}) is the $L^1(A)$ limit of $\vp_\eps$, so for $\L^1$-a.e.\ $x \in A_\nu$, $\vp_x$ coincide with the restriction of $\vp$ to $A_x$. Therefore, since the vector $\nu$ is taken arbitrarily, $\vp \in BV(A)$ (see Proposition 6.9 in \cite{ABS}), and since $A$ is any open relatively compact subset of $\D$, we derive that $\vp \in BV_{loc}(\D)$. Using (\ref{slicingineq}), (\ref{lw1d3x}) Fatou's lemma and Fubini's formula, we also obtain

 \ba
 	\liminf_{\eps \to 0} \eps \F_\eps (v_\eps,\vp_\eps) &\geq& \int_{A_\nu}  \F(1,\varphi_x \, ;A_x) \, dx \nonumber \\
				&=&    \frac{2\sigma}{\pi} \, \int_{A_\nu}  dx\, \int_{A_x} \rho_x^{\nicefrac32} d |D\vp| \nonumber \\
				&=&   \frac{2\sigma}{\pi} \, \int_\D \rho_x^{\nicefrac32} d (\L^1 \mesrest A_\nu \otimes |D\vp_x| \mesrest A_x) \label{estimatenu} \,.
 \ea

 Now, for every $\eps>0$, let $\mu_\eps$ be the energy distribution in $\D$ associated with the pair $(v_\eps,\vp_\eps)$, that is, the positive Radon measure which for every Borel set $E \subset \R^2$ is given by

\ben
	\mu_\eps(E) =  \eps \F_\eps (v_\eps,\vp_\eps \,; E \cap \D) \,.
\een

 From (\ref{unifboundF}),  the total mass $\|\mu_\eps\|$ is uniformly bounded. De La Vallée Poussin compactness criterion (see \cite{AFP}, page 26) gives then that (up to a subsequence)  $\mu_\eps$ converges weakly$^*$  to some finite measure $\mu$ on $\D$. We claim that

\ben
	\mu\geq  2\sigma \rho^{\nicefrac32} \cdot \H^1 \mesrest S\vp \,.
\een

We will prove this using Besicovitch derivation Theorem (see \cite{AFP}, page 54). First, after (\ref{unifboundF}) for every $K \subset\subset\D$ there is $R_K \in (0,\lambda)$ such that

\be \label{0infty}
	0 \leq \mu(K) \leq \mu(B(0,R_K)) \leq \liminf_{\eps\to0} \mu_\eps(B(0,R_K))  < \infty.
\ee

Hence,  $\mu$ is a positive Radon measure in $\D$, and for $\H^1$-a.e.\ $x \in S\vp$ the limit

\be \label{deff}
	f(x) = \lim_{r\to0^+} \frac{\mu(B_r(x))}{\H^1(B_r(x) \cap S\vp)}
\ee

exists, and we have

\be \label{propf}
	\mu \geq f \cdot \H^1 \mesrest S\vp \,.
\ee

 Let $x_0 \in S\vp \cap A$.   Since $A \subset\subset\D$, $\overline{B_r(x_0)} \subset\D$ for $r$ small enough. We assume\footnote{In fact this holds for all $r$ except countably many (see \cite{AFP}, page 29).} moreover that $\mu(\partial B(x_0,r))=0$. Proposition 1.62 in \cite{AFP} and estimate  (\ref{estimatenu}) yield

\ba
	\mu(B_r(x_0)) &=& \lim_{\eps \to 0}\mu_\eps(B_r(x_0)) \nonumber \\
	&\geq&  \frac{2\sigma}{\pi} \, \int_{B_r(x_0)} \rho_x^{\nicefrac32} d (\L^1 \mesrest A_\nu \otimes |D\vp_x| \mesrest A_x)\nonumber \\
	&\geq& \frac{2\sigma}{\pi} \, \inf_{B_r(x_0)} \rho_x^{\nicefrac32}  \,  \int_{B_r(x)} d (\L^1 \mesrest A_\nu \otimes |D\vp_x| \mesrest A_x) \,. \label{mu(Br(x)}
\ea

In Proposition \ref{lowerboundonlines} we proved that for $\L^1$-a.e.\ $x \in A_\nu$, $\vp_x \in SBV(A_x) \cap L^\infty(A)$, and

\ben
	\int_{A_\nu} dx\int_{A_x} \H^0(S_{\vp_x}) < \infty \,.
\een

Hence, after Theorem 2.3 in \cite{AT}, $\vp \in SBV(A) \cap L^\infty(A)$ and

 \be \label{equalityB0}
	\int_{B_r(x_0)} d (\L^1 \mesrest A_\nu \otimes |D\vp_x| \mesrest A_x) =  \pi \int_{B_r(x_0)} |\langle \nu_\vp , \nu \rangle  | \, d \H^1 \,,
\ee

where $\nu_\vp$ is the measure theoretic inner normal to the Caccioppoli set $\{ \vp =\pi\}$. Putting (\ref{equalityB0}) in (\ref{mu(Br(x)}) we obtain

\ban
	\mu(B_r(x_0)) &\geq&   2\sigma  \, \inf_{B_r(x_0)} \rho_x^{\nicefrac32}  \,  \int_{B_r(x_0)} |\langle \nu_\vp , \nu \rangle  | \, d \H^1  \\
			     &\geq&    2\sigma   \, \inf_{B_r(x_0)} \rho_x^{\nicefrac32}  \,      \inf_{B_r(x_0) \cap S\vp} |\langle \nu_\vp , \nu \rangle  | \,  \H^1(B_r(x_0)\cap S\vp)  \,.
\ean

Since $S\vp$ is a rectifiable set in $A$, for $\H^1$-a.e.\ $x_0 \in S\vp$, $\nu_\vp(x)$ is continuous in $B_r(x_0)$ for $r$ small enough. Thus, taking $\nu=\nu_\vp(x_0)$ and since $\rho$ is continuous, we get

\ben
	\lim_{r\to0^+} \frac{\mu(B_r(x_0))}{\H^1(B_r(x_0) \cap S\vp)} \geq  2\sigma   \rho^{\nicefrac32}(x_0) \,
\een

for $\H^1$-a.e.\ $x_0 \in S\vp$. Hence, (\ref{deff}) and (\ref{propf}) yield the claim. The definition of $\mu$ gives then

 \ben
 	\liminf_{\eps \to 0} \eps \F_\eps (v_\eps,\vp_\eps \,; A) = \liminf_{\eps \to 0}  \mu_\eps(A) \geq \mu(A)  \geq 2\sigma \, \int_{S\vp \cap A} \rho^{\nicefrac32}  d\H^1  \,.
 \een

Finally, taking an increasing sequence $\{A_k\}_{k \in \N} $ with $A_k \subset\subset \D$, we get

 \ben
 	\liminf_{\eps \to 0} \eps \F_\eps (v_\eps,\vp_\eps) \geq 2\sigma    \int_{S\vp \cap \D } \rho^{\nicefrac32}  d\H^1 \,,
 \een

which gives (\ref{lowerboundF}). \\ \qed \\

%%%%%%%%%
%%%%%%%%%%
\subsection{Remark about a lower bound for $v_\eps$ in the transition zone}

We end this section with a discussion about the infinimum of $v_\eps$ in the transition zone. Let $\{(v_\eps,\vp_\eps) \}_{\eps >0}$ be a sequence of minimizers of $ \F_\eps $, and let $\vp \in BV_{\text{loc}}(\D \,;\, \{ 0, \pi \})$ be the $L^1_{\text{loc}}$-limit of $\vp_\eps$ given in (\ref{compactnessgcvtheo}). Let $K \subset \subset \D $ be an open smooth set, with non negligible intersection with $S\vp$, that is,

\ben
	\H^1(K \cap S\vp) > 0 \,.
\een

For every $\eps >0$, we define

\ben
	m_{\eps,K} =   \inf_{x \in K}  v_\eps(x) \,.
\een

We would like to obtain an upper bound for $m_{\eps,K}$, in connection with an open question in \cite{BeLinWeiZhao}, namely
\be\label{upmeps}
	m_{\eps,K}  \leq C_K  (g_\eps \eps^2)^{-\nicefrac14} \,.
\ee

If we assume that we have the upper and lower inequalities for each $\eps>0$, that is

\be \label{finiteremarklimsup}
	 \eps F_\eps( \tilde v_\eps) \leq \F(\vp) \,
\ee

and

\be \label{finiteremarkliminf}
	 \eps F_\eps(\tilde v_\eps) \geq \F(\vp) \,,
\ee

we can give estimates on $G_\eps$ in order to obtain the upper bound for $m_{\eps,K}$. So assume that we have (\ref{finiteremarklimsup}) and (\ref{finiteremarkliminf}). On the one hand, estimates (\ref{estetaK}) and  (\ref{estrhopD}) give then

\ban
	 \eps G_\eps(v_\eps, \vp_\eps \,;\, K) &\geq& \frac14 \sqrt{g_\eps \eps^2} \, m_{\eps,K}^3 \left( \inf_K \rho^{\nicefrac32} -  C_K \eps^2 |\ln \eps^2| \right) \int_{K} |\nabla \vp_\eps| \sin \vp_\eps  \\
	 &\geq& C_K \, \sqrt{g_\eps \eps^2} \, m_{\eps,K}^3 \int_{K} |\nabla \vp_\eps| \sin \vp_\eps  \,.
\ean

We claim that the integral here below is bounded away from zero. Indeed, if this not the case, we will have

\ben
	\liminf_{\eps \to 0} \int_{K} |\nabla \vp_\eps| \sin \vp_\eps  = 0\,.
\een

Hence, since $\vp_\eps \to \vp$ in $L^1(K)$,  the coarea formula together with the lower semi continuity of the $BV$ norm imply the contradiction

\ben
	0 =  \liminf_{\eps \to 0}  \int_0^\pi \sin t \, dt \int_{K} |D \textbf{1}_{ \{ \vp_\eps <t \} }|  \geq  \int_0^\pi \sin t \, dt    \int_{K}   |D \textbf{1}_{ \{\vp = 0 \}}| = 2\,  \H^1( S\vp \cap K) \,.
\een

We thus derive that there is $C'_K >0$ such that

\be  \label{remarkB}
	 \eps G_\eps(v_\eps, \vp_\eps \,;\, K) \geq C'_K \, \sqrt{g_\eps \eps^2} \, m_{\eps,K}^3  \,.
\ee

In the other hand, by inspection of the proof of Proposition \ref{propUB} (see estimate (\ref{limepsGeps})), we see that the pair of test function $(\tilde v_\eps, \tilde \vpeps)$ satisfies

\be  \label{remarkD}
	\eps G_\eps(\tilde v_\eps, \tilde \vp_\eps  \,;\, K) \leq C  (g_\eps \eps^2)^{-\nicefrac14}  \,.
\ee

Hence, considering (\ref{finiteremarklimsup})-(\ref{remarkD}), together with the fact that $(v_\eps, \vp_\eps )$  minimizes $\F_\eps$, we obtain

\ban
	 2\sigma \int_{K} \rho^{\nicefrac32} |D\vp| + C_K \, \sqrt{g_\eps \eps^2} \, m_{\eps,K}^3  &\leq&   \eps \F_\eps( v_\eps,   \vp_\eps \,;\, K) \\
	  &\leq&   \eps \F_\eps( \tilde v_\eps,  \tilde \vp_\eps \,;\, K) \\
	  &\leq&  2\sigma \int_{K} \rho^{\nicefrac32} |D\vp|  + C  (g_\eps \eps^2)^{-\nicefrac14} \,.
\ean

Multiplying both sides of the previous inequality by $(g_\eps \eps^2)^{\nicefrac14}$ we find the  upper bound (\ref{upmeps}) for $m_{\eps,K}^3 $.

However, we are not able to prove  (\ref{finiteremarklimsup}) and (\ref{finiteremarkliminf}) as such because of the error terms.
Indeed, the proof of the upper bound of Theorem \ref{theo} says that there is a sequence $\{( \tilde v_\eps, \tilde\vp_\eps) \}_{\eps >0}$ such that

\be \label{remarklimsup}
	\limsup_{\eps \to 0} \eps F_\eps( \tilde v_\eps) \leq \F(\vp) \,.
\ee

In the proof of (\ref{remarklimsup}), we first approximate the locally Cacciopoli set $A=\{ \vp = \pi\}$ by characteristics functions of open sets $A_k$ with compact smooth boundary. This gives a small error in terms of $k \in \N$ in the upper bound inequality (\ref{remarklimsup}). Then for each $k \in \N$, we construct a test function for which (\ref{remarklimsup}) holds, up to a small error term depending on a parameter $\delta>0$.  In these two steps, we use diagonal extraction arguments in order to get rid of the error terms, so it is not possible to compute them explicitly. Similarly, in the proof of the lower bound of Theorem \ref{theo}, we use the compactness of bounded Radon measures, so we cannot estimate the error term in the lower bound inequality

\be \label{remarkliminf}
	\liminf_{\eps \to 0} \eps F_\eps(\tilde v_\eps) \geq \F(\vp) \,.
\ee

%%%%%%%%%%%%%%%%
%%%%%%%%%%%%%%%%
\section{Proof of the $\Gamma$-convergence for $\eps \left( \E_\eps(\cdot) - E_\eps(\eta_\eps) \right)$} \label{gammacv}

%%%%%%
\subsection{Proof of the compactness and the lower bound inequality in Theorem \ref{theo}:}

let $\{(u_{1,\eps}, u_{2,\eps})\}_{\eps>0}$ be a sequence of minimizers of $\E_\eps$ in $\mathcal{H}$ satisfying (\ref{upgcvtheo}). From Proposition \ref{propu1u2vphi}\textbf{(i)}, the pairs $(v_\eps, \vp_\eps)$ are well defined by (\ref{chvariables}), belong to $ Lip_{loc}( \R^2 \,;\, (0,+\infty) \times [0,\pi])$ and satisfy (\ref{boundlinftv2}). Proposition \ref{energysplittingprop} yields $\eps \F_\eps(v_\eps,\vp_\eps)< \infty$. Thus, the hypotheses of Proposition \ref{propLB} are fulfilled by $(v_\eps, \vp_\eps)$ and we have

 \ban
	(v_\eps,\vp_\eps) \to (1,\vp) \quad &\text{ in }& \quad L^1_{loc}(\D) \times L^1_{loc}(\D)
\ean

with $\vp \in X$, and

\ben
	\liminf_{\eps\to0} \eps \F_\eps(v_\eps,\vp_\eps\, ) \geq \F(\varphi \, ) \,.
\een

Equality (\ref{splitener}) yields then

\ben
	\liminf_{\eps\to0} \eps \left( \E_\eps(u_{1,\eps}, u_{2,\eps}) - E_\eps(\eta_\eps) \right) \geq \F(\varphi ) \,.
\een

Finally, using identity (\ref{invert}) we get

\ben
	(u_{1,\eps},u_{2,\eps}) \to \sqrt{\rho} \, ( \mathbf{1}_{\{ \vp =0\}}, \mathbf{1}_{\{ \vp = \pi\}} ) \quad \text{in} \quad L^1_{loc}(\D) \times  L^1_{loc}(\D)  \,.
\een

\qed \\

In order to prove the upper bound we have to work a little more. We first modify the pairs of test functions from Proposition \ref{propUB} to make them satisfy the mass constraints (\ref{mass2}). We prove then that this modification do not change the limit of the energy. We finish by verifying that the pairs of modified test functions are the image by (\ref{chvariables}) of a pair in $\H$, and we conclude using Proposition \ref{propUB}.  \\

%%%%%%
\subsection{Proof of the upper bound inequality in Theorem \ref{theo}:}

\textit{(Step 1 : Modification of the pairs of test functions)} With the notations from the proof of Proposition \ref{propUB}, we write $N_{\eps} = N_{\eps R_\delta}$ and we define $(\check{v}_\eps, \vpeps)$ the sequence of pairs of test functions such that

\be  \label{mathcalFeps2}
	\limsup_{\eps \to 0} \eps \F_\eps (\check{v}_\eps, \vpeps) \leq \F(\vp)  \,.
\ee

Consider $\kappa \in C^\infty(\R_+ \,;\, [0,1])$ with $supp \, \kappa \subset (0,1)$ and $\kappa =1$ in $(0,\nicefrac12)$. Since $A$ is a non empty open set,  there is  $B_0 = B_{r_0}(x_0) \subset \subset A\cap \D$.  For $\ell \in[-1,1]$ and  $\tau \in (\nicefrac12,1)$, define $\kappa_\eps=\kappa_{\eps,\ell,\tau}$ by

\ben
	\kappa_\eps(x) =  \eps^\tau \ell \,  \kappa(\nicefrac{|x-x_0|}{r_0}) \,.
\een

We define then $ \hat{v}_\eps = \check{v}_\eps + \kappa_\eps$ and $ v_\eps = c_\eps \hat{v}_\eps$, with $c_\eps = \| \eta_\eps \hat{v}_\eps \|_2^{-2}$. For $\eps$ small enough $N_\eps $ and  $B_0$ are disjoints. We estimate

\ban
	c_\eps^{-1} &=& 1+ \int_{N_\eps \cup B_0} \eta_\eps^2 ( \hat{v}^2_\eps-1) \\
			   &=& 1 + 2 \int_{B_0} \eta_\eps^2 \kappa_\eps \, + \int_{N_\eps} \eta_\eps^2 (  \check{v}^2_\eps-1)  \, + \,  \int_{B_0} \eta_\eps^2 \kappa^2_\eps \\
			   &=& 1+ 2 \int_{B_0} \eta_\eps^2 \kappa_\eps \, + \O(\eps)  + \O(\eps^{2\tau}) \,.
\ean

Hence, using that $\tau \in (\nicefrac12,1)$ we get $c_\eps^2 = 1 - r_\eps$ with

\be \label{estreps}
	r_\eps = 4 \int_{B_0} \eta_\eps^2 \kappa_\eps \, + \O(\eps) = \O(\eps^\tau) \,.
\ee

Notice that for $\eps$ small enough, $r_\eps $ may be positive or negative depending on the sign of $\ell$. \\

The definition of $ w^y_{\eps,T}$ insures that $v_\eps>0$. The first mass constraint in (\ref{mass2}) is immediately satisfied by the definition of $v_\eps$. Remember the definition of $\vp_\eps$ in (\ref{defphieps}). For the second mass constraint we write

\ban
	c_\eps^{-2} \int_{\R^2} \eta^2_\eps v^2_\eps \cos \vp_\eps &=&  \int_{\R^2} \eta^2_\eps ( \mathbf{1}_{\R^2 \backslash(A \cup N_\eps)} - \mathbf{1}_{ A \backslash (N_\eps \cup B_0)} + \mathbf{1}_{ N_\eps \cup B_0} \, \hat{v}^2_\eps \cos \vp_\eps) \,.
\ean

Adding and removing $ \mathbf{1}_{ N_\eps \backslash A} \, \eta^2_\eps $, $ \mathbf{1}_{ N_\eps \cup A} \, \eta^2_\eps $ and  $ \mathbf{1}_{ B_0} \, \eta^2_\eps $ in the previous integral, we get

\ba \label{int1}
	c_\eps^{-2} \int_{\R^2} \eta^2_\eps v^2_\eps \cos \vp_\eps &=&   \int_{\R^2} \eta^2_\eps ( \mathbf{1}_{\R^2} -  2\mathbf{1}_{ A} )  +   \int_{B_0}  \eta^2_\eps (\check{v}_\eps + \kappa_\eps)^2-1  \\
	&& + \int_{N_\eps}  \eta^2_\eps ( \check{v}_\eps^2  \cos \vp_\eps   -\mathbf{1}_{ A} +\mathbf{1}_{\R^2 \backslash A}) \,.
\ea

For the third term in (\ref{int1}), we have that $\eta_\eps$, $\check{v}_\eps$ and  $\cos \vp_\eps$ are bounded while $\mathcal{L}^2(N_\eps) = \O(\eps)$. Hence,

\ba \label{int2}
	\int_{N_\eps}  \eta^2_\eps ( \check{v}_\eps^2  \cos \vp_\eps   -\mathbf{1}_{ A} +\mathbf{1}_{\R^2 \backslash A})  = \O(\eps) \,.
\ea

For the first term in (\ref{int1}), using that $\int_{\R^2}  \eta^2_\eps = 1 = \alpha_1+\alpha_2$ and that $\int_{\D \cap A}\rho= \alpha_2 $, we obtain

\ban
	\int_{\R^2} \eta^2_\eps ( \mathbf{1}_{\R^2} -  2\mathbf{1}_{ A} ) = \alpha_1-\alpha_2 + \int_{A \cap \D} (\eta_\eps^2 - \rho) + \int_{A \backslash \D} \eta_\eps^2   \,.
\ean

Using (\ref{estetanearboundary}) we get, for $\alpha \in (\nicefrac12,\nicefrac35)$ and $\gamma \in (\nicefrac12,\nicefrac34)$,

\ba
	 \int_{A \cap \D} (\eta_\eps^2 - \rho) &=& \int_{A \cap B(0,\lambda-\eps^\alpha)}  (\eta_\eps^2 - \rho) +  \int_{(A \cap \D) \backslash B(0,\lambda-\eps^\alpha)}  (\eta_\eps^2 - \rho) \nonumber  \\ \nonumber  \\
	 &=& \O(\eps^\gamma)  + \O(\eps^\alpha) \,.   \label{int3}
\ea

Moreover, from (\ref{estetarad}), we have $\eta^2_\eps(x) \leq \eta^2_\eps(x_\alpha)$ in $A \backslash \D$, with $x_\alpha \in \partial B(0,\lambda-\eps^\alpha)$. From  (\ref{estetanearboundary}) and (\ref{estrhopD}) we get

\ben
	\eta^2_\eps(x) \leq \eta^2_\eps(x_\alpha) = \eta_\eps^2(x_\alpha)  - \rho(x_\alpha)  + \rho(x_\alpha) =  \O(\eps^\alpha) \,,
\een

so using that $A$ is a bounded set we obtain

\be  \label{int4}
	 \int_{A \backslash \D}  \eta_\eps^2  = \O(\eps^\alpha) \,.
\ee

For the second term in (\ref{int1}), the definitions of $\kappa_\eps$ and $r_\eps$ yield

\ba \label{int5}
	\int_{B_0}  \eta^2_\eps \kappa_\eps (2 + \kappa_\eps ) = \frac12 r_\eps + \O(\eps) + \O(\eps^{2\tau}) \,.
\ea

Putting (\ref{int2})-(\ref{int5}) in (\ref{int1}) and considering (\ref{estreps}) we get

\ben
	c_\eps^{-2} \int_{\R^2} \eta^2_\eps v^2_\eps \cos \vp_\eps = \alpha_1 - \alpha_2 +   \frac12 r_\eps + \O(\eps^\beta)  \,,
\een

where $\beta = \min\{1,\alpha,\gamma,2\tau \} =  \min\{\alpha,\gamma\} \in (\nicefrac12,\nicefrac35) $. Hence, (\ref{estreps}) gives

\ben
	\int_{\R^2} \eta^2_\eps v^2_\eps \cos \vp_\eps   - (\alpha_1 - \alpha_2) = \left(\frac12 -  (\alpha_1 - \alpha_2)   \right)  r_\eps + \O(\eps^\beta) \,.
\een

Suppose now, without loss of generality, that $\alpha_1 - \alpha_2 \leq \nicefrac12$. The definition of $r_\eps$ and $\kappa_\eps$, together with (\ref{estetaK}), (\ref{estrhopD})  and $B_0  \subset \subset A\cap \D$, give then

\ban
	|r_\eps| &\geq& 4 \inf_{B_0} \eta^2_\eps  \int_{B_0} \kappa_\eps^2  + \O(\eps)   \\
	 	&\geq& 4 \inf_{B_0}  \eta^2_\eps \,  |\ell|  \,  \eps^\tau  \,  \int_{B_{\nicefrac {r_0}2}(x_0)} \kappa_\eps^2+ \O(\eps) \\
		 & \geq& c  \,  |\ell|  \,  \eps^\tau  + \O(\eps)  \,,
\ean

for some $c>0$ not depending on $\eps$. Hence, if we take $\ell=1$ in the definition of $\kappa_\eps$,  for $\eps$ small enough we have

\ben
	\int_{\R^2} \eta^2_\eps v^2_\eps \cos \vp_\eps   - (\alpha_1 - \alpha_2) \geq c'  \eps^\tau ( 1 + \eps^{1-\tau} -  \eps^{\beta-\tau})  \,.
\een

Analogously,  taking now $\ell=-1$, we get

\ben
	\int_{\R^2} \eta^2_\eps v^2_\eps \cos \vp_\eps   - (\alpha_1 - \alpha_2) \leq c^{''}  \eps^\tau ( -1 + \eps^{1-\tau}+  \eps^{\beta-\tau})  \,.
\een

Since $\beta \in (\nicefrac12,\nicefrac35)  $, we can choose $\tau \in (\nicefrac12,\beta)$ and obtain

\ben
	\int_{\R^2} \eta^2_\eps v^2_\eps \cos \vp_\eps   >  \alpha_1 - \alpha_2  \quad \text{ if } \quad \ell =1
\een

 and

\ben
	\int_{\R^2} \eta^2_\eps v^2_\eps \cos \vp_\eps   <  \alpha_1 - \alpha_2    \quad \text{ if }  \quad \ell =-1  \,.
\een

Hence, there exists $\ell_\eps \in (-1,1)$ such that  for $\eps$ small enough, the associated pair $(v_\eps , \vp_\eps )$ satisfy the second mass constraint in (\ref{mass2}). \\

%%%%%%%%%%%%
\textit{(Step 2 : Computing the energy)}. We now compute the energy of $(v_\eps , \vp_\eps )$. We recall that $N_{\eps \tilde t_\eps}$ is the transition zone of $\vp_\eps$ defined in (\ref{defphieps}). For the energy $G_\eps$, we have that $\vp_\eps$ is constant out of $N_{\eps \tilde t_\eps}$, while $ v_\eps = c_\eps \check{v}_\eps$ in $N_{\eps \tilde t_\eps}$ with $c_\eps = 1+\O(\eps^\tau)$. Hence,

\be \label{Geps2}
	\eps G_\eps (v_\eps, \vpeps)  = (1+\O(\eps^\tau)) \eps G_\eps (\check{v}_\eps, \vpeps)   \,.
\ee

For the energy $F_\eps$, we have that $v_\eps = c_\eps (1+\kappa_\eps)$ in $B_0$. The definition of $\kappa_\eps$ gives then, $|\nabla v_\eps|^2 = \O(\eps^{2\tau}) $ and $\{ 1 - v_\eps^2\}^2 = \O(\eps^{2\tau})$. Hence,

\be \label{intB0}
	\eps F_\eps (v_\eps \,;\, B_0)  = \O(\eps^{2\tau-1}) = o_{\eps \to 0}(1) \,.
\ee

In $\R^2\backslash(N_\eps \cap B_0)$, we have that $v_\eps = c_\eps$, so $|\nabla v_\eps| =0 $ and $\{ 1 - v_\eps^2\}^2 = \O(\eps^{2\tau})$. As before we get

\be \label{intVepsB0}
	\eps F_\eps (v_\eps \,;\, \R^2\backslash(N_\eps \cap B_0))  = \O(\eps^{2\tau-1}) = o_{\eps \to 0}(1) \,.
\ee

In $N_\eps$, we have that $v_\eps = c_\eps \check{v}_\eps$. Hence, $|\nabla v_\eps|^2 = (1+\O(\eps^{\tau})) |\nabla  \check{v}_\eps|^2 $ and $\{ 1 - v_\eps^2\}^2 =(1+\O(\eps^{\tau}))   \{ 1 - \check{v}_\eps^2\}^2 +\O(\eps^{\tau})$, which gives

\ba
	\eps F_\eps (v_\eps \,;\,  N_\eps)  &=& (1+\O(\eps^{\tau})) \,  \eps F_\eps (\check{v}_\eps \,;\,  N_\eps) + \O(\eps^{\tau}) \eps^{-1} \L^2(N_\eps)  \nonumber \\
	&=& (1+\O(\eps^{\tau})) \,  \eps F_\eps (\check{v}_\eps \,;\,  N_\eps) + o_{\eps \to 0}(1)  \label{intVeps} \,.
\ea

Since $\check{v}_\eps$ is constant out of $N_\eps$, we have $F_\eps (\check{v}_\eps) = F_\eps (\check{v}_\eps \,;\,  N_\eps) $. Putting together (\ref{mathcalFeps2}) and (\ref{Geps2})-(\ref{intVeps}), we obtain

\be \label{limsupagain}
	\limsup_{\eps \to 0} \eps \F_\eps (v_\eps, \vpeps) = \limsup_{\eps \to 0} \eps \F_\eps (\check{v}_\eps, \vpeps) \leq \F(\vp)  \,.
\ee

%%%%%%%%
\textit{(Step 3 : identification of $(v_\eps,\vp_\eps)$)} The pairs of test functions satisfies the hypothesis from Proposition \ref{propu1u2vphi}\textbf{(ii)}, so defining $(u_{1,\eps},u_{2,\eps})$  by (\ref{invert}) we have $(u_{1,\eps},u_{2,\eps}) \in \mathcal{H}$ and $u_{1,\eps}^2+u_{2,\eps}^2 >0$. Hence, after Proposition (\ref{energysplittingprop}) relation (\ref{splitener}) holds, and (\ref{limsupagain}) yield

\ben
	\limsup_{\eps \to 0} \eps \left( \E_\eps(u_{1,\eps}, u_{2,\eps}) - E_\eps(\eta_\eps) \right)  =\limsup_{\eps \to 0} \eps \F_\eps (v_\eps, \vpeps)  \leq \F(\vp)  \,.
\een

\qed \\

%%%%%%%%%
\subsection{Proof of Corollary \ref{coro}:} Let $\tilde \vp \in X$ with $\F(\tilde\vp) < +\infty$. From the upper bound inequality in Theorem \ref{theo}, there is a sequence $(\tilde u_{1,\eps}, \tilde u_{2,\eps}) \in \mathcal{H}$ such that

\ben
	 \limsup_{\eps \to 0} \eps \left( \E_\eps( \tilde u_{1,\eps}, \tilde u_{2,\eps}) - E_\eps(\eta_\eps) \right) \leq \F( \tilde \vp)  \,.
\een

Since $(u_{1,\eps}, u_{2,\eps}) $ minimize $\E_\eps$ in $\mathcal H$, the previous inequality yields

\be \label{prooftheo1a}
	 \limsup_{\eps \to 0} \eps \left( \E_\eps(   u_{1,\eps},   u_{2,\eps}) - E_\eps(\eta_\eps) \right) \leq \F( \tilde \vp)  \,,
\ee

so in particular  $(u_{1,\eps}, u_{2,\eps}) $  satisfy (\ref{upgcvtheo}). Hence, from the compactness and the lower bound inequality in Theorem \ref{theo}, there is $\varphi \in X$ and a subsequence $(u_{1,\eps'}, u_{2,\eps'}) $ with

\ben
	\liminf_{ \eps' \to 0} \eps' \left( \E_{\eps'} (u_{1,\eps'}, u_{2,\eps'}) - E_{\eps'}(\eta_{\eps'})\right) \geq \F(\vp) \,.
\een

This inequality is verified for every subsequence of $( u_{1,\eps},  u_{2,\eps})$, so we have

\be \label{prooftheo1b}
	\liminf_{\eps \to 0} \eps \left( \E_\eps(   u_{1,\eps},   u_{2,\eps}) - E_\eps(\eta_\eps) \right) \geq   \F(\vp)  \,.
\ee

From (\ref{prooftheo1a}) and (\ref{prooftheo1b}), we obtain

\be \label{prooftheo1c}
	\F(\tilde \vp) \geq \limsup_{\eps \to 0} \eps \E_{\eps} (    u_{1,\eps},    u_{2,\eps}) \geq  \liminf_{\eps \to 0} \eps \E_{\eps} ( u_{1,\eps},  u_{2,\eps}) \geq \F(\vp) \,,
\ee

so $\F(\vp) = \inf_{X} \F$. Taking $\tilde \vp = \vp$ in (\ref{prooftheo1c}) yields

\ben
	\lim_{\eps \to 0} \eps\left( \E_{\eps} (  u_{1,\eps},  u_{2,\eps}) - E_\eps(\eta_\eps) \right) =  \inf_{X} \F  \,.
\een

 \qed \\

%%%%%%%%%%%%%%%%%%%%%%%
 \subsection{Proof of Corollary \ref{corobreak}  }

We start  proving that when $\alpha_1 $ is not to close to $0$ or $ 1$, the minimizers of $\F$ in $X$ are not radially symmetric. We show that for any radially symmetric $\vp \in X$, $\F(\vp ) > \F(\vp_{ds}) $, where the support of $\vp_{ds} \in X$ is a disk sector. We first prove this for functions such that $\{\vp=0\}$ is a disk or an annulus. Then, we generalize by induction the result  to radial functions such that $\{\vp=0\}$ is composed of a finite number of connected components.  We conclude then by approximating any radially symmetric  $ \vp \in X$ by this kind of functions.\\

 We recall that $\rho$ is given in (\ref{rho})  and that $X$ is the space of functions $ \vp \in BV_{loc}(\D;\{0,\pi\})$ such that

\be \label{rapelmc}
	\int_{ \{\vp =0\} } \rho = \alpha_1 \,.
 \ee

If $\vp_{ds} \in X$ is such that $\{\vp_{ds}=0\}$ is a disk sector, we easily compute

\ben
	\frac{\F(\vp_{ds})}{8\sigma} =   \frac3{16}\,.
\een

 For $0\leq R^-\leq R^+ \leq \lambda$ we denote $A(R^-,R^+)$ the annulus of center the origin, inner radius $R^-$ and outer radius $R^+$. \\

 If $\vp_\alpha \in X$ is such that $\{\vp=0\}=A( 0,R_\alpha )$ and $\int_{ A(0, R_\alpha )} \rho = \alpha$, then $R_\alpha = \lambda(1-\sqrt{ \alpha})^{\nicefrac12}$ and

 \be  \label{Fvpalphaf}
	\frac{\F(\vp_\alpha)}{8\sigma} =   f(\alpha)	\,,
\ee

where $f : [0,1] \to \R_+$ is the concave function $f(\alpha) =(1- \alpha )^{\nicefrac34} (1-\sqrt{1-\alpha })^{\nicefrac12}$. We see that there exists

\ben
	  \delta_0  \approx  0.1486
\een

	such that if $\alpha  \in [\delta_0,1-\delta_0]$, then $f(\alpha )>\nicefrac3{16}$.

 \begin{prop} \label{symbreakF}
If $\alpha_1 \in [\delta_0,1-\delta_0]$, then the minimizers of $\F$ in $X$ are not radially symmetric.
\end{prop}

\textbf{Proof:}  \textit{(Step 1)}  Let $R \in (0,\lambda)$ and consider  $   \vp_1^d \in X$ such that $\{ \vp_1^d  = 0\} = A(0,R)$. From (\ref{rapelmc}),  we have that $\F(\vp_1^d  ) / 8 \sigma= f(\alpha_1)$ so (\ref{Fvpalphaf})  yields

\be \label{energyvp0}
 	\F(\vp_1^d  ) >\F(\vp_{ds}) \,.
\ee

Since $\alpha_2=1-\alpha_1 \in [\delta_0,1-\delta_0]$, the similar inequality holds if $\{\vp_1^d = 0\} = A( R ,\lambda)$. \\

Consider now $   \vp_1^a \in X$ such that $\{ \vp_1^a  = 0\} = A(R_1,R_2)$, with $0<R_1<R_2<\lambda$. Writing

\ben
	\beta_1 = \int_{A(0,R_1)} \rho  \,, \quad \quad  \beta_2 = \int_{A( R_1,R_2)} \rho \quad \text{ and } \quad \beta_3 = \int_{A(  R_2,\lambda)} \rho  \,,
 \een

 we compute

\ben
	\frac{\F(\vp_1^a)}{8\sigma} =   f(\beta_1)	+ f(\beta_1+\beta_2)\,.
\een

From (\ref{rapelmc}), we have that  $\beta_2=\alpha_1$ and $\beta_1+\beta_3 = \alpha_2$  so

\ben
	\frac{\F(\vp_1^a)}{8\sigma} =   f(\beta_1)	+ f(\beta_1+\alpha_1)	\,.
\een

The right hand size of the previous equality is a concave function of $\beta_1$ and the value of $\beta_1$ may vary  between $0$ and $\alpha_2$. If $\beta_1=0$ then   $ \F(\vp_1^a)  / 8\sigma=  f(\alpha_1)$. If $\beta_1= \alpha_2$, since $\alpha_1+\alpha_2=1 $ we find $ \F(\vp_1^a)  / 8\sigma=   f(\alpha_2)$. We derive

\be  \label{energyvp01}
	\F(\vp_1^a)  >\F(\vp_{ds}) \,.
\ee

  \textit{(Step 2)} Let $n \in \N^*$ and consider $\vp_n \in X$ such that

\ben
	\{\vp_n = 0\} = \bigcup_{j=1}^n A_{2j }    \,,
\een

with $A_{2j } = A( R_{2j } ^- ,R_{2j } ^+ )$ and

\ben
	 0\leq R_{2j-2}^- < R_{2j-2}^+ < R_{2j }^-  <R_{2j  }^+\leq \lambda
\een

for $2\leq j\leq n$. We write $\beta_{2j} = \int_{A_{2j}} \rho$, $\beta_1= \int_{A(0,R_2^-)} \rho$, $\beta_{2n+1}= \int_{A( R_{2n}^+,\lambda)} \rho$  and

\ben
	\beta_{2j+1} = \int_{ A( R_{2j } ^+ ,R_{2j+2 } ^- )} \rho
\een

for $ 1\leq j \leq n-1$. Notice that we allow $A(0,R_2^-)$ or $A( R_{2n}^+,\lambda)$  to be empty, but this only implies that  $\beta_1=0$  or $\beta_{2n+1}=0$. With this notation, we have

\be \label{hypind1}
	 \sum_{i=1}^{n}    \beta_{2i} = \alpha_1 \,, \quad    \quad    \sum_{i=1}^{n}    \beta_{2i+1} = \alpha_2
 \ee

and

\ben
	\frac{\F(\vp_n)}{8\sigma}=   \sum_{j=1}^{2n } f\left( \sum_{i=1}^{j} \beta_i  \right) =: g_n(\beta_1,\cdots, \beta_{2n}) \,.
\een

By induction, we are going to prove the following property:

\ban
	(\mathcal P_n)  \quad && \forall \, \beta_1, \cdots ,\beta_{2n+1} \in [0,1]  \text{ such that  }   \sum_{i=1}^{n}    \beta_{2i} =  \alpha_1 \text{ and }    \sum_{i=1}^{n}    \beta_{2i+1} =  \alpha_2 \,,  \\  &&  g_n(\beta_1,\cdots, \beta_{2n}) > \frac{\F(\vp_{ds})}{8\sigma}\,.
\ean	

If $n=1$ we are in one of the three cases  analyzed in Step 1, so (\ref{energyvp0}) and (\ref{energyvp01}) yield $(\mathcal P_1)$. \\

 Let us assume that $(\mathcal P_n)$ holds and consider $ \beta_1, \cdots ,\beta_{2n+3} \in [0,1]$ such that

 \be \label{2nplus3}
 	\sum_{i=1}^{ n+1}  \beta_{2i} = \alpha_1 \quad \text{ and } \quad  \sum_{i=1}^{ n+1}  \beta_{2i+1} = \alpha_2 \,.
\ee

We have

\ben
	g_{n+1}( \beta_1,\cdots, \beta_{2n+2}) =  \sum_{j=1}^{2n  } f\left( \sum_{i=1}^{j} \beta_i  \right) + f\left( \sum_{i=1}^{ 2n+1 } \beta_{ i} \right)   + f\left( \sum_{i=1}^{ 2n+2 } \beta_{ i} \right)   \,.
\een

The right hand side of the previous equality is a concave function of $\beta_{2n+2}$. The value of $\beta_{2n+2}$ may vary between $0$ and $\alpha_1$.  Suppose first that $\beta_{2n+2}=0$.  Then, defining

\ben
	\tilde \beta_j  =  \beta_j  \quad  \text{ if }   \quad  j=1,\cdots,2n \qquad \text { and } \qquad  \tilde \beta_{2n+1}   = \beta_{2n+1}+ \beta_{2n+3}  \,,
 \een

the $\tilde \beta_i$'s satisfy (\ref{hypind1}) and  we have

\ben
	g_{n+1}( \beta_1,\cdots, \beta_{2n+2}) \geq  \sum_{j=1}^{2n  } f\left( \sum_{i=1}^{j}   \beta_i  \right)  = g_n(\tilde \beta_1,\cdots, \tilde \beta_{2n}) \,.
\een

 Hence, $(\mathcal P_n)$ yields $  g_{n+1}( \beta_1,\cdots, \beta_{2n+2}) >  \F(\vp_{ds})  / 8\sigma$. \\

 Suppose now  that  $\beta_{2n+2}= \alpha_1$.   From (\ref{2nplus3}) this implies $ \beta_{2j}=0  $ for every $j=1,\cdots, n$. Then, defining

\ben
	\tilde \beta_1  = \sum_{j=1}^{2n+1} \beta_j \,,    \qquad  \tilde \beta_2 =\beta_{2n+2}  \qquad \text { and } \qquad \tilde \beta_3=  \beta_{2n+3}  \,,
 \een

the $\tilde \beta_i$'s satisfy (\ref{hypind1}) and  we have

\ban
	g_{n+1}( \beta_1,\cdots, \beta_{2n+2}) &\geq&    f\left( \sum_{i=1}^{ 2n+1 } \beta_{ i} \right)   + f\left( \sum_{i=1}^{ 2n+2 } \beta_{ i} \right) \\
	&=& f(\tilde \beta_1) + f(\tilde \beta_1 +\tilde \beta_2) \\
	&=&  g_1(\tilde \beta_1,\tilde \beta_2)   \,.
\ean

Hence, $(\mathcal P_1)$ yields $ g_{n+1}( \beta_1,\cdots, \beta_{2n+2}) >  \F(\vp_{ds})/8\sigma $. We derive that the result holds for all the possible values of $\beta_{2n+2}$. \\

 We have proved that if $\vp_n \in X$ is radial and its support has a finite number of connected components, then

\be \label{inequaln}
	\F(\vp_n)  >  \F(\vp_{ds}) \,.
\ee

  \textit{(Step 3)}  Suppose now that $\vp \in X$ is a radially symmetric function such that $\{ \vp = 0\}$ has an infinite number of connected components. Since $\vp$ has locally finite perimeter in $\D$, $\{ \vp = 0\}$  is the union of a countable family of disjoints annuli. We write

\ben
	 \{ \vp = 0\} = \bigcup_{j \in \mathbb{Z}}  A_{2j }
\een

 with  $A_{2j } = A(R_{2j }^- ,R_{2j  }^+)  $ such that

 \be \label{A2j}
	 0< R_{2j }^- < R_{2j }^+ < R_{2j+2 }^-  <R_{2j +2 }^+< \lambda  \,.
\ee

 For every $n \in \N$, we define a function $\vp_n : \D \to \{ 0,\pi\}$ by

 \ben
 	 \{ \vp_n = 0 \} = 	\bigcup_{j  =-n}^n  A_{2j } \bigcup \tilde A_{2n+2} \bigcup \tilde A_{-2n-2} \,,
 \een

such that

\ben
	\tilde A_{2n+2} = A(L_n^+,\lambda) \quad \text{ and } \quad  \tilde A_{-2n-2} = A( 0,L_n^- )
 \een

with  $ L_n^-,L_n^+ $ to be chosen next. If $(L_n^-,L_n^+) =   (0,\lambda)$, then (\ref{A2j}) gives

\ben
	\int_{\{\vp_n=0\}} \rho =\sum_{ j  =-n}^n \int_{ A_{2j }} \rho   <   \int_{\{\vp =0\}} \rho  \,.
\een

Similarly if $(L_n^-,L_n^+) = (R_{-2n}^-,R_{2n}^+)$, then

\ben
	\int_{\{\vp_n=0\}} \rho =	 \sum_{ j  \in \mathbb{Z} } \int_{ A_{2j }} \rho  +   \sum_{ j\geq n} \int_{A(  R_{2j}^+, R_{2j+2}^- ) }  \rho +\sum_{ j\leq - n } \int_{A(  R_{2j-2}^+, R_{2j }^- )}  \rho    >  \int_{\{\vp =0\}} \rho \,.
\een

Hence, by continuity there is a pair $(L_n^-,L_n^+) \in(0,R_{-2n}^-) \times(R_{2n}^+, \lambda) $ such that $\int_{\{\vp_n=0\}} \rho = \int_{\{\vp =0\}} \rho =\alpha_1 $. Clearly $\vp_n \in BV_{loc}(\D)$, so $\vp_n \in X$. Moreover, (\ref{A2j})  yields

\be \label{limLplusLminus}
	\lim_{n \to \infty} L_n^- = 0 \quad \text{ and } \quad \lim_{n \to \infty} L_n^+ = \lambda \,.
\ee

We have

\ben
	\F(\vp)= \sum_{j \in \mathbb Z}  \int_{ \partial B(0, R_{2j}^+)} \rho^{\nicefrac32}  \, d\H^1\,,
 \een

and  since $\rho$ is radially symmetric

\ban
	\F(\vp_n) &=& \sum_{j=-n}^n  \int_{ \partial B(0, R_{2j}^+)} \rho^{\nicefrac32}  \, d\H^1  +  2 \pi \left(  \rho^{\nicefrac32} (L_n^+) \, L_n^ +  +   \rho^{\nicefrac32} (L_n^-)  \, L_n^- \right)  \,.
\ean

 From (\ref{limLplusLminus}), the  last term in the previous equality goes to zero as $n \to +\infty$, so $\lim_{n \to \infty} 	 \F(\vp_n)=  \	\F(\vp ). $ Hence, since $\{\vp_n=0\}$ has a finite number of  connected components, (\ref{inequaln})  yields $\F (\vp) >   \F ( \vp_{ds}) $, which ends the proof. \\ \qed  \\

\textbf{Proof of Corollary \ref{corobreak}:} Suppose that   $\a_1 \in [\delta_0,1-\delta_0]$ and that $\{ (u_{1,\eps}, u_{2,\eps}) \}_{\eps>0} $ is a sequence of radially symmetric pairs such that $(u_{1,\eps}, u_{2,\eps})$ minimizes	$ \E_\eps$  under the mass constraints (\ref{mass}). Then, $\varphi_\eps$ defined by  (\ref{relationvpu12})   is also radially symmetric. Consider $\varphi_{\eps,0}$, the restriction of $\varphi_\eps$ to a slice of $\D$ passing through $0$. From Proposition \ref{lowerboundonlines}, $\varphi_{\eps,0}$ belongs to $SBV_{loc}([ 0,\lambda]  \,;\, \{0,\pi \})$ and converges in  $L^1_{loc}([0,\lambda])$  to  $\varphi_{0}$. Hence, $\varphi_\eps$ converges  in  $L^1_{loc}(\D)$  to the radial function $\varphi$ given by $\varphi (x) = \varphi_{0}(|x|) $. From Corollary \ref{coro}, we know that $\varphi$ minimizes $\F$ over $X$, which yields a contradiction with Proposition \ref{symbreakF}. \\ \qed \\

%%%%%%%%%
\section{Appendix}

We end this article given the proof of Lemma \ref{approxA}, which is essentially the same of Lemma 4.3 in \cite{BOU}, which in turn is a generalization of Lemma 1 in \cite{Mod}. For completeness we give here the details of the proof.\\

\textbf{Proof of Lemma \ref{approxA}:} \textit{(Step 1)} Suppose first that $\D \cap A$ and $\D \backslash A $ have both non empty interior and let

\be \label{B12}
	\overline{B(x_1,\delta)} \subset \D \cap A \qquad \text{ and } \qquad \overline{B(x_2,\delta)} \subset \D \backslash A \,.
\ee

We first approximate $A$ by sets of finite perimeter in $\D$. For $k\geq2$ we define $\D_k = \D \cap B(0,\lambda(1-\nicefrac1k))$ and $A'_k = A \cap \D_k$. We have that $\partial_* A'_k \subset (\partial_* A \cap \D_k) \cup \partial \D_k$, so

\ben
	 \int_{\partial_* A'_k} \rho^{\nicefrac32} \, d\H^1 \leq  \int_{\partial_* A \cap \D_k} \rho^{\nicefrac32} \, d\H^1 +  \int_{\partial \D_k} \rho^{\nicefrac32} \, d\H^1 \,.
\een

Using Lebesgue dominated convergence theorem, the first term in the right hand side of the inequality converges to  $\int_{\partial_* A} \rho^{\nicefrac32} \, d\H^1 < + \infty$. The definition of $\D_k $ and (\ref{estrhopD}) yield

\ban
  \int_{\partial \D_k} \rho^{\nicefrac32} \, d\H^1 \leq \| \rho\|_{L^\infty(\partial \D_k)}^{\nicefrac32} \, \H^1(\partial \D_k) \leq   \Big(\frac {2\lambda^2} k \Big)^{\nicefrac32} \, \H^1(\partial \D) = o_{k \to \infty}(1) \,.
\ean

Hence,

\be \label{limenergyAprime}
	\lim_{k \to \infty} \int_{\partial_* A'_k} \rho^{\nicefrac32} \, d\H^1 \leq \int_{\partial_* A} \rho^{\nicefrac32} \, d\H^1 \,.
\ee

 %%%
\textit{(Step 2)} Since $A'_k$ has finite perimeter in $\D$, it can be approximated (see the proof of Lemma 1 in \cite{Mod}) by open bounded sets $\tilde A_k$, such that

\ba
	 \L^2(\tilde A_k \Delta A'_k) &\leq& \frac1k  \label{AkL2} \\
	A'_k \subset \tilde A_k  + B(0,\nicefrac1k)   &&  \text{ and}  \qquad \tilde A_k \subset  A'_k  + B(0,\nicefrac1k)  \label{AktildeB}\\
	\H^1(\partial \tilde A_k \cap \pD) &=& 0 \label{H1tildeAk} \,.
\ea

The definition of $A'_k$ and (\ref{AkL2}) imply \textbf{(i)}. Using (\ref{B12}) and (\ref{AktildeB}), for large enough $k$ we get

\be \label{Bx1Bx2}
	B(x_1,\delta) \subset \tilde A_k \qquad \text{ and } \qquad B(x_2,\delta) \subset \D \backslash  \tilde A_k \,.
\ee

Moreover, using (ii) from Proposition 2.3 in \cite{BOU} and the fact that $\tilde A_k$ belongs to a sequence $\tilde A^n_k$ such that $\|\tilde A^n_k\|_{BV(\D)} \to \|A'_k\|_{BV(\D)} $ as $n \to 0$, we have

\be \label{energytildeAkprime1k}
	 \int_{\D} \rho^{\nicefrac32} |D \mathbf{1}_{\tilde A_k}| \leq  \int_{\D} \rho^{\nicefrac32} |D \mathbf{1}_{ A'_k}| + \frac1k \,.
\ee

 Also, the definition of $A'_k$ and (\ref{AkL2}) yield

 \be \label{gammak}
 	\gamma_k := \int_{\tilde A_k} \rho -  \int_{A} \rho =o_{k \to \infty}(1) \,.
 \ee

 %%%
\textit{(Step 3)}  Now, we set

\ben
	A_k = \left\{ \begin{array}{ccc}
		\tilde A_k \backslash B(x_1,r_{1,k}) & \text{ if } & \gamma_k> 0 \\
		\tilde A_k   & \text{ if } & \gamma_k = 0 \\
		\tilde A_k \cup B(x_1,r_{2,k}) & \text{ if } & \gamma_k<0
	\end{array}\right. \,,
\een

where $r_{1,k}$ and $r_{2,k}$ are chosen to satisfy

\ben
	\int_{B(x_1,r_{1,k})} \rho = \int_{B(x_2,r_{2,k})} \rho  = 	\gamma_k \,.
\een

Since $r \mapsto  \int_{B(x_{1,2},r)} \rho$ is continuous and decreasing for $r \in (0, \delta)$, $r_{1,k}$ and $r_{2,k}$ are unique and tend to zero as $k \to \infty$. Then, we derive from  (\ref{Bx1Bx2}) and (\ref{gammak}), for large enough $k$, that

 \ben
	 \int_{\D \cap A_k} \rho =  \int_{\tilde A_k} \rho  - \gamma_k=  \int_{A} \rho \,.
 \een

Moreover, from (\ref{H1tildeAk}) and (\ref{Bx1Bx2}), we have $\H^1(\partial A_k) = 0$ for $k$ large enough, so \textbf{(ii)} is proved. Using again (\ref{Bx1Bx2}) we obtain

\be \label{ineqA1AB12}
	 \int_{\D} \rho^{\nicefrac32} |D \mathbf{1}_{ A_k}| \leq  \int_{\D} \rho^{\nicefrac32} |D \mathbf{1}_{ \tilde A_k}| + \| \rho \|_\infty \, H^1(\partial B(x_1,r_{1,k})  \cup \partial B(x_2,r_{1,2}) ) \,.
\ee

Hence, using (\ref{energytildeAkprime1k}), we obtain

\ben
	 \int_{\D} \rho^{\nicefrac32} |D \mathbf{1}_{A_k}| \leq  \int_{\D} \rho^{\nicefrac32} |D \mathbf{1}_{ A'_k}| + o_{k \to \infty} (1) \,,
\een

so (\ref{limenergyAprime}) gives

\ben
	 \limsup_{k \to \infty} \int_{\D} \rho^{\nicefrac32} |D \mathbf{1}_{A_k}| \leq  \int_{\D} \rho^{\nicefrac32} |D \mathbf{1}_{ A}| \,
\een

We have proved \textbf{(iii)}.\\

%%%
\textit{(Step 4)} We now remove the condition that $\D \cap A$ and $\D \backslash A $ have no empty interior. First, we notice that $\L^2( \D \cap A) =0$ and $\L^2( \D \backslash A) =0$ are not possible because of the mass constraints in (\ref{mass2}). Hence, there exists $x_1$ a point of density of $\D \cap A$ and $x_2$ a point of density of $\D \backslash A$. Consider the function

\ben
	\Phi(\delta_1,\delta_2) = \int_{A_{12}} \rho - \int_{A} \rho  \,,
\een

where $A_{12} = A \cup B(x_1,\delta_1) \backslash B(x_2,\delta_2)  $. Since $\rho > 0$ in $\D$, for any $\delta >0$ we have

\ben
	\Phi(\delta,0) > 0 \qquad \text{ and } \qquad  \Phi(0,\delta) < 0 \,.
\een

Since $\Phi$ is continuous, there is $t=t_\delta \in (0,1)$  such that $\Phi(t\delta,(1-t)\delta)=0$. Define $A_\delta =  A \cup B(x_1,(1-t)\delta) \backslash B(x_2,t\delta) $ and $\vp=\pi \mathbf{1}_{A_\delta}$. Both $\D \cap A_\delta$ and $\D \backslash A_\delta $ have no empty interior and $\int_{A_\delta} \rho = \int_{A} \rho$. Moreover $\L^2(A_\delta  \Delta A) \to 0$ as $\delta \to 0$, and using an inequality similar to (\ref{ineqA1AB12}), we get

\ben
	\limsup_{\delta \to 0}   \int_{\D} \rho^{\nicefrac32} |D \mathbf{1}_{A_\delta}| \leq   \int_{\D} \rho^{\nicefrac32} |D \mathbf{1}_{A}| \,.
\een

Finally, for each $A_\delta$ we apply the construction from steps 1-3 and conclude thanks a diagonal argument, see Corollary 1.16 in \cite{attouch}. \\ \qed \\

%%%%%%%%%%%%%%%%
\noindent{\bf Acknowledgements} The second author would like to acknowledge discussions with Guy Bouchitté, Pierre Seppecher and Duvan Henao. We would like to thank Clément Gallo.

\bibliographystyle{acm}
\bibliography{bibaajrl}

\begin{thebibliography}{10}

\bibitem{Aflivre}
{\sc Aftalion, A.}
\newblock {\em Vortices in Bose-Einstein Condensates}, vol.~67 of {\em Progress
  in Nonlinear Differential Equations and Their Applications}.
\newblock Birkh\"auser, 2006.

\bibitem{AJR}
{\sc Aftalion, A., Jerrard, R.~L., and Royo-Letelier, J.}
\newblock Non-existence of vortices in the small density region of a
  condensate.
\newblock {\em J. Funct. Anal. 260\/} (2011), 2387--2406.

\bibitem{alb}
{\sc Alberti, G.}
\newblock Variational models for phase transitions, an approach via
  {$\Gamma$}-convergence.
\newblock {\em Calculus of Variations and Differential Equations, Springer,
  Berlin, 2000\/}, 95--114.

\bibitem{ABS}
{\sc Alberti, G., Bouchitt{\'e}, G., and Seppecher, P.}
\newblock Phase transition with the line-tension effect.
\newblock {\em Arch. Rational Mech. Anal. 144}, 1 (1998), 1--46.

\bibitem{AFP}
{\sc Ambrosio, L., Fusco, N., and Pallara, D.}
\newblock {\em Functions of bounded variation and free discontinuity problems}.
\newblock Oxford New York : Clarendon Press, 2000.

\bibitem{AT}
{\sc Ambrosio, L., and Tortorelli, V.~M.}
\newblock Approximation of functionals depending on jumps by elliptic
  functionals via {$\Gamma$}-convergence.
\newblock {\em Comm. Pure Appl. Math. 43}, 8 (1990), 999--1036.

\bibitem{attouch}
{\sc Attouch, H.}
\newblock {\em Variational convergence for functions and operators}.
\newblock Pitman Advanced Publishing Program, 1984.

\bibitem{BeLinWeiZhao}
{\sc Berestycki, H., Lin, T.-C., Wei, J., and Zhao, C.}
\newblock On phase-separation model: Asymptotics and qualitative properties.
\newblock {\em Arch. Rational Mech. Anal.\/} (2013), to appear.

\bibitem{BeTer}
{\sc Berestycki, H., Terracini, S., Wang, K., and Wei, J.}
\newblock On entire solutions of an elliptic system modeling phase separations.
\newblock {\em Preprint\/} (2012), to appear.

\bibitem{BOU}
{\sc Bouchitt\'e, G.}
\newblock Singular perturbations of variational problems arising from a
  two-phase transition model.
\newblock {\em Appl. Math. Optim. 21}, 3 (1990), 289--314.

\bibitem{braides}
{\sc Braides, A.}
\newblock {\em Approximation of free-discontinuity problems}.
\newblock Lecture Notes in Mathematics, Vol. 1694. Springer, 1998.

\bibitem{brezis1}
{\sc Brezis, H.}
\newblock Semilinear equations in {${\bf R}^N$} without condition at infinity.
\newblock {\em Appl. Math. Optim. 12}, 3 (1984), 271--282.

\bibitem{CaffLin2}
{\sc Caffarelli, L.~A., and Lin, F.-H.}
\newblock Singularly perturbed elliptic systems and multi-valued harmonic
  functions with free boundaries.
\newblock {\em J. Amer. Math. Soc. 21}, 3 (2008), 847--862.

\bibitem{ctv3}
{\sc Conti, M., Terracini, S., and Verzini, G.}
\newblock On a class of optimal partition problem related to the
  {F}u$\check{\text{c}}$\'ik spectrum and to the monotonicity formulae.
\newblock {\em Calc. Var. Partial Differential Equations 22}, 1 (2005), 45--72.

\bibitem{EG}
{\sc Evans, L.~C., and Gariepy, R.~F.}
\newblock {\em Measure Theory and Fine Properties of Functions}.
\newblock CRC Press, 1992.

\bibitem{Ga}
{\sc Gallo, C.}
\newblock Expansion of the energy of the ground state of the
  {G}ross--{P}itaevskii equation in the {T}homas--{F}ermi limit.
\newblock {\em ArXiv e-prints\/} (May 2012).

\bibitem{GaPe}
{\sc Gallo, C., and Pelinovsky, D.}
\newblock On the {T}homas--{F}ermi ground state in a harmonic potential.
\newblock {\em Asymptotic Analysis 73\/} (2011), 53--96.

\bibitem{Giusti}
{\sc Giusti, E.}
\newblock {\em Minimal Surfaces and Functions of Bounded Variation}.
\newblock Monographs in Mathematics. Birkh{\"a}user Boston, 1984.

\bibitem{hall}
{\sc Hall, D., Matthews, M., Wieman, C., and Cornell, E.}
\newblock Measurements of relative phase in binary mixtures of
  {B}ose-{E}instein condensates.
\newblock {\em Phys. Rev. Lett. 81\/} (1998), 1543--1547.

\bibitem{IM}
{\sc Ignat, R., and Millot, V.}
\newblock The critical velocity for vortex existence in a two-dimensional
  rotating {B}ose-{E}instein condensate.
\newblock {\em J. Funct. Anal. 233\/} (2006), 260--306.

\bibitem{KaSou}
{\sc {Karali}, G.~D., and {Sourdis}, C.}
\newblock {The ground state of a Gross-Pitaevskii energy with general potential
  in the Thomas-Fermi limit}.
\newblock {\em ArXiv e-prints\/} (May 2012).

\bibitem{LM}
{\sc Lassoued, L., and Mironescu, P.}
\newblock Ginzburg-{L}andau type energy with discontinuous constraint.
\newblock {\em J. Anal. Math. 77\/} (1999), 1--26.

\bibitem{DM}
{\sc Maso, G.~D.}
\newblock Integral representation on {BV}($\omega$) of {$\Gamma$-limits} of
  variational integrals.
\newblock {\em Manuscripta Mathematica 30}, 4 (1979), 387--416.

\bibitem{MaAf}
{\sc Mason, P., and Aftalion, A.}
\newblock Classification of the ground states and topological defects in a
  rotating two-component {B}ose-{E}instein condensate.
\newblock {\em Phys. Rev. A 84}, 3 (2011), 033611.

\bibitem{McC}
{\sc McCarron, D.~J., Cho, H.~W., Jenkin, D.~L., K\"oppinger, M.~P., and
  Cornish, S.~L.}
\newblock Dual-species {B}ose-{E}instein condensate of $^{87}\mathrm{Rb}$ and
  $^{133}\mathrm{Cs}$.
\newblock {\em Phys. Rev. A 84\/} (2011), 011603.

\bibitem{Mod}
{\sc Modica, L.}
\newblock The gradient theory of phase transitions and the minimal interface
  criterion.
\newblock {\em Arch. Rational Mech. Anal. 98}, 2 (1987), 123--142.

\bibitem{modugno}
{\sc Modugno, G., Modugno, M., Riboli, F., Roati, G., and Inguscio, M.}
\newblock A two atomic species superfluid.
\newblock {\em Phys. Rev. Lett. 89\/} (2002), 190404--190408.

\bibitem{NoTaTeVe}
{\sc Noris, B., Tavares, H., Terracini, S., and Verzini, G.}
\newblock Uniform {H}{\"o}lder bounds for nonlinear {S}chr\"odinger systems
  with strong competition.
\newblock {\em Comm. Pure Appl. Math. 63}, 3 (2010), 267--302.

\bibitem{JRL}
{\sc Royo-Letelier, J.}
\newblock Segregation and symmetry breaking of strongly coupled two-component
  {B}ose-{E}instein condensates in a harmonic trap.
\newblock {\em Calc. Var. Partial Differential Equations\/} (2012), to appear.

\bibitem{WeWe1}
{\sc Wei, J., and Weth, T.}
\newblock Asymptotic behaviour of solutions of planar elliptic systems with
  strong competition.
\newblock {\em Nonlinearity 21}, 2 (2008), 305--317.

\end{thebibliography}

\end{document}